\documentclass[a4paper]{article}
\usepackage[english]{babel}
\usepackage[utf8x]{inputenc}
\usepackage{amsmath}
\usepackage{graphicx}
\usepackage{mathtools}
\usepackage{microtype}
\usepackage{mathtools}
\usepackage{breqn}
\usepackage{empheq}
\usepackage{tensor}
\usepackage{array}
\usepackage{multirow}
\usepackage{transparent}
\usepackage{fontenc}
\usepackage{booktabs}
\usepackage[numbers]{natbib}
\usepackage{amssymb}
\usepackage{soul}
\usepackage{xcolor}
\usepackage{bibentry}
\usepackage{algorithm}
\usepackage{algorithmic}
\usepackage{amsmath}
\usepackage{amsfonts}
\usepackage{array}
\usepackage{subfig}
\usepackage{booktabs} 
\usepackage{geometry} 
\usepackage{amssymb}
\usepackage[colorinlistoftodos]{todonotes}
\usepackage{hyperref}
\usepackage{cleveref}
\usepackage{geometry} 
\usepackage{lipsum}
\allowdisplaybreaks
\numberwithin{equation}{section}
\newtheorem{theorem}{Theorem}[section]
\newtheorem{corollary}{Corollary}[theorem]

\newtheorem{lemma}[theorem]{Lemma}
\newtheorem{example}[theorem]{Example}

\newtheorem{Remark}{Remark}[theorem]
\usepackage{cleveref}
\begin{document}
\title{Uniform norm error estimate for rectangular finite element approximation of a 2D turning point problem}
\author{
Shallu$^{*}$, Sudipto Chowdhury$^{*}$, Vikas Gupta$^{*}$ \\[2pt]
$^{*}$\tiny{Centre for Mathematical and Financial Computing, The LNM Institute of Information Technology, Jaipur-302031, Rajasthan, India}
}
\maketitle
\textbf{Abstract}
\\This work presents error analysis for a finite element method applied to a two-dimensional singularly
perturbed convection-diffusion turning point problem. Utilizing a layer-adapted Shishkin mesh, we prove uniform convergence in the maximum norm in the x-layer regions and $\varepsilon$-independent bounds for the coarse region. The analysis, critically based on the properties of a discrete Green's function, guarantees the method's robustness and accuracy in capturing sharp solution layers.
\section{INTRODUCTION}\label{section1}
This paper examines a particular class of convection-diffusion boundary value problems characterized by homogeneous Dirichlet boundary conditions. The general form of the problem under consideration is
\begin{subequations}\label{problem}
    \begin{align}
    \ -\varepsilon \Delta u +  \textbf{b}(x,y).\nabla u + c(x,y) u & = f(x,y), \quad  \text{in} \quad \Omega = (-1,1) \times (-1,1), \\
    u& = 0\quad \text{on} \quad \Gamma = \partial \Omega.
\end{align} 
\end{subequations}
Here, $\varepsilon$ is a small positive perturbation parameter  $0<\varepsilon<<1$. The coefficients $b =(b_1,b_2)$, $ c$, and the source term $f(x,y)$ are assumed to be sufficiently smooth over the closure of the domain $\bar{\Omega}=\Omega \cup \partial \Omega$.\\
A particular case of the general problem \eqref{problem} is of special interest and constitutes the focus of the present study. In particular, the convection vector $b(x,y)$ is specified such that $b_1(x,y) = xa(x,y)$ and $b_2(x,y)=0 \; \forall \; (x,y) \in \Omega$. This specification leads to the following simplified problem: 
 \begin{subequations}\label{current problem}
      \begin{align}
             -\varepsilon \Delta u + x a(x,y) u_x + c(x,y) u &= f(x,y), \quad \text{in} \quad\Omega = (-1,1) \times (-1,1),\\
              u &= 0,\quad \text{on}  \quad \Gamma = \partial \Omega.\\
              |a(x,y)| \geq \alpha > 0, \quad |c(x,y)|& \geq 0, \quad  c - \frac{1}{2} \frac{\partial b_1}{\partial x} \geq \gamma >0.
        \end{align}
 \end{subequations}
  In the resulting problem \eqref{current problem}, the convection coefficient vanishes along the line $x=0$; therefore, the problem is referred to as a turning-point problem. Here, we make a very crucial assumption about $\alpha$, namely that $\alpha\geq5/2$. Consequently, $\alpha$ cannot be chosen arbitrarily small. For \eqref{current problem}, an interior layer occurs at $x=0$. The solution may not exhibit boundary layers at the lateral boundaries ($x=-1\; \text{and} \;x=1$); instead, it displays parabolic-type layers along the top and bottom boundaries ($y=1$) and ($y=-1$).\\
 This paper focuses on the finite element analysis of a two-dimensional singularly perturbed problem with turning points along the line $x=0$, i.e., a turning line defined on the domain (-1,1)$\times$(-1,1). The combination of a turning point in a two or higher-dimensional setting presents a significant challenge for both theoretical analysis and numerical approximation. To address this, we employ a conforming quadrilateral finite element on a properly constructed Shishkin mesh that is adapted to the layer structures induced by the perturbation parameter and the turning point.\\
Researchers have been actively involved in the progression of finite difference schemes for singularly perturbed turning point problems (SPTPPs). Their work presents both single and multiple interior turning points; for example, see \cite{Becher_25, Berger_84, Farrell_88, Gupta_21, Geng_13, Kadalbajoo_10, Kumar_19, Natesan_98, Natesan_03}. Berger et al.\cite{Berger_84} and Farrell \cite{Farrell_88} conducted analyses of SPTPPs involving a single interior turning point. Their work confirmed a priori bounds on the analytic solution and its derivatives, which were crucial for demonstrating the uniform convergence of their respective numerical scheme. Specifically, Berger et al.\cite{Berger_84} derived an $\varepsilon$-uniform error estimate for a modified El-Mistikawy-Werle finite difference scheme. In parallel, Farrell \cite{Farrell_88} achieved uniform convergence estimates for a broad category of upwinding-type schemes. In a unique contribution, Natesan and Ramanujam \cite{Natesan_98} introduced a composite numerical approach for SPTPPs that exhibit twin boundary layers. Their method combines a classical numerical scheme with an exponential fitted-difference scheme to generate an approximate solution. To implement this, the authors divided the problem domain into four distinct sub-regions and applied an asymptotic expansion approximation to solve the problem within each segment. Kadalbajoo and Gupta \cite{Gupta_21} contributed by deriving asymptotic bounds for the derivatives of the continuous solution to SPTPPs featuring twin boundary layers. These bounds were subsequently used to establish second-order $\varepsilon$-uniform convergence for a B-spline collocation method applied to a piecewise uniform layer-adapted mesh. Natesan et al.\cite{ Natesan_03} proved almost first-order 
$\varepsilon$-uniform convergence for their proposed scheme, which was based on the classical upwind finite difference method (FDMs) implemented on a Shishkin mesh. Separately, Geng and Qian \cite{Geng_13} introduced a stretching variable method and a reproducing kernel method for addressing SPTPPs whose solutions present twin boundary layers. In the work by Becher and Roos \cite{Becher_25}, the Richardson extrapolation method was applied to an upwind scheme, utilizing a piecewise layer-adapted mesh. This application successfully refined the order of accuracy from $N^{-2}\ln N$ to $N^{-2}(\ln N)^2$, where $N$ represents the number of intervals in the spatial direction.\\
Compared to the abundant literature on finite difference methods (FDMs), research on applying the finite element method (FEMs) to SPTPPs is limited. However, FEMs generally provide a more efficient route to higher-order convergence and require less stringent data regularity than FDMs.\\
Over the past decades, researchers have increasingly directed their efforts toward developing FEMs for SPTPPs. A key milestone was achieved in 1994 when Sun and Stynes \cite{Sun94} presented and rigorously analyzed piecewise linear Galerkin FEMs. This scheme, tested on diverse piecewise equidistant meshes, was designed for SPPs with interior turning points and solutions featuring interior layers. They successfully established parameter-uniform convergence for their method in the standard $L^2$-norm and a weighted energy norm. Building upon this work, Becher \cite{Becher_16} in 2016 introduced a higher-order finite element approximation for SPPs with interior turning points. Employing Liseikin's \cite{Liseikin18} mesh, Becher's contribution was significant in obtaining error bounds that were independent of the perturbation parameter in the energy norm. Work by Becher \cite{Becher18} focused on a one-dimensional SPTPPs that exhibited interior layer characteristics in its solution. The problem was discretized via a standard Galerkin FEMs, then stabilized with SDFEM on piecewise equidistant meshes, yielding parameter-independent error bounds in both the energy and SD-norms. Subsequently, in \cite{Aasna2024, Aasna Rai2, Chen08} singularly perturbed boundary turning point problems, both with and without delay, were analyzed. The application of SDFEM in these studies achieved parameter-independent convergence, approaching second order, in the maximum norm. Aasna and Rai \cite{Aasna2024} further developed the SDFEM framework for time-dependent SPPs with boundary turning points. Their strategy involved employing an $\theta$-scheme on an equidistant time mesh and SDFEM on a spatial Shishkin mesh, successfully establishing parameter-uniform stability and convergence estimates in the maximum norm for their method. Ranjan \cite{Rajeev_Ranjan} considered (\ref{current problem}), and a uniform energy norm estimate of order $N^{-1}\log^2 N$ is derived for the SIPG method on a layer-adapted Shishkin mesh.
\\
To the best of our knowledge, the pointwise error estimate of the finite element method for two-or higher-dimensional SPTPPs has received considerably less attention in the existing literature.  In \cite{Stynes}, Stynes used a Galerkin finite element method on a piecewise equidistant mesh to solve a singularly perturbed boundary value problem of convection diffusion type in two dimensions without a turning point. That analysis established uniform convergence with respect to the perturbation parameter of order $N^{-1} \log N$ in the global energy norm, and of order $N^{-1/2}\log^{3/2} N$ pointwise near the outflow boundary, but not in the coarse region. A key strength of the analysis lies in the use of a discrete Green’s function, which enables the derivation of pointwise error estimates. This observation strongly motivates the present work. In \cite{zhang}, pointwise error estimates based on discrete Green's functions for the SDFEM were established by Zhang for the convection diffusion problems containing only boundary layers. Though the analysis was based on a uniform mesh. Pointwise energy norm and $L^2$-norm estimates on layer-adapted meshes for purely boundary layer and turning point convection diffusion problems were derived in \cite{lins} for SDFEM. There, the barrier function technique is used instead of the discrete Green's function technique. 
\\\\
In this work, we investigate pointwise error estimates for the classical finite element method applied to problems of type (\ref{current problem}) on a Shishkin mesh. We employ the standard finite element method on a piecewise uniform (Shishkin) mesh, incorporating specific transition parameters along the $x$- and $y$-directions, given respectively by 
\[
\lambda_x = \min\left(\frac{2\varepsilon}{\alpha} \log\left(\frac{1}{\varepsilon}\right),\frac{1}{2}\right), \quad \lambda_y = \min\left(2\sqrt{\frac{\varepsilon}{\beta}} \log\left(\frac{1}{\varepsilon^{3/2}}\right),  \frac{1}{4}\right).
\]
This framework is used to solve convection-diffusion boundary value problems with homogeneous Dirichlet boundary conditions and a turning point in the domain $(-1,1)\times(-1,1)$. We establish pointwise uniform convergence with respect to the perturbation parameter, achieving an order of $\frac{(\log(1/\varepsilon))^{1/2}}{N^{1/2}}$ in the $x$-layer region. To the best of our knowledge, this is the first work to provide a pointwise norm error estimate for turning point problems employing the classical finite element method on the layer-adapted meshes.  \\\\
This paper proceeds as follows: In Section \ref{section1}, we establish bounds for the solution of the continuous problem (\ref{current problem}). In Section \ref{section2}, we present the classical finite element method, describing its key components such as the finite element space, the chosen basis functions, and the layer-adapted piecewise uniform mesh (known as the Shishkin mesh) with specific mesh transitions. Due to the symmetry along the $y$-axis only need to analyze the problem in the domain $[0,1]\times[-1,1]$ along with we decompose our domain into four parts $\Omega_{c}, \;\Omega_{f,x}, \;\Omega_{f,y}$ and $\Omega_{f,xy}$ represents coarse region, $x$-layer region, $y$-layer region and $xy$-layer region respectively. In Section \ref{section3}, we begin by recalling the discrete Green's function, which plays a key role in deriving pointwise error estimates. Then, we proceed to derive $L^2$-norm estimates for the discrete Green's function. A key contribution of this section is a $\varepsilon$-independent estimate for the discrete Green's function when the source point lies in either the coarse region or the $x$-layer region. This estimate serves as a critical component in the proof of our main theorem. In Section \ref{section4}, we establish $L^{\infty}$-error (pointwise) estimates for the diffusion term. To achieve this, we decompose the solution into four parts: R, $F_{1}$, $F_{2}$, and $F_{12}$ corresponding to no-layer (smooth) region, $x$-layer, $y$-layer, and $xy$-layer regions, respectively. Then, we recall an interpolation theorem from \cite{Rajeev_Ranjan}, followed by the assumption that source points lie in coarse and $x$-layer regions, and we derived pointwise error estimates for all layer regions. In section \ref{section5}, following the interpolation framework by Zarin and H. Ross, we develop $L^\infty$-norm error (pointwise) estimates for the convection and reaction term. The analysis yields $\varepsilon$-independent error bounds provided the source point resides within the coarse and x-layer regions. In section \ref{section6}, combining estimates for convection, diffusion, and reaction components, we establish the main theorem, which confirms uniform convergence in the perturbation parameter ($\varepsilon$) specifically within these regions. Based on the estimates for the convection, diffusion, and reaction terms, we derive the main theorem of this work, which demonstrates uniform convergence in the perturbation parameter $\frac{(\log(1/\varepsilon))^{1/2}}{N^{1/2}}$  pointwise specifically $(x_m,y_n)\in\Omega_{f,x}$ and $\varepsilon$-independent bounds when $(x_m,y_n)\in\Omega_{c}$. In \cref{section-6}, based on the theoretical, we present an example of the convection diffusion turning point problem. The double mesh principle is employed to compute the error estimates and determine the rate of convergence for various values of the perturbation parameter $\varepsilon$ and the mesh discretization parameter N. The pointwise error estimates are presented in Table \ref{tab:table1} and Table \ref{tab:table2} when source points are within smooth and $x$-layer regions, respectively. Table \ref{tab:table3} provides the convergence behavior in both subdomains. The numerical solution for Example \ref{example} is depicted in the accompanying figure, thereby demonstrating the effectiveness of the proposed approach. To the best of our knowledge, there has not been any work on the pointwise error estimate for the two or higher-dimensional turning point problems.  
\subsection{The continuous problem}
The weak formulation for  \eqref{current problem} is to find an unique $u \in H_0^1(\Omega) = \{u\in H^1(\Omega) ; u|_{\partial \Omega}=0\}$ such that 
\begin{equation*}
    \varepsilon \int_\Omega \nabla u. \nabla v \; dx dy  + \int_\Omega x a(x,y) u_x v \; dxdy  + \int_\Omega c(x,y) u v \; dxdy = \int_\Omega f(x,y) v \; dxdy \quad \forall \; v\in  H_0^1(\Omega)
\end{equation*}
 \begin{equation} \label{weak formul}
             B(u,v) = (f,v)  \quad \forall \; v,u \in H^1_0(\Omega),
         \end{equation}
    where $B(.,.)$ is defined as 
         \begin{equation} \label{weak}
             B(u,v) = (\varepsilon \nabla u, \nabla v) + (\textbf{b}.\nabla u , v) +(cu,v) \quad \forall \;  u,v \in H^1_0(\Omega).
         \end{equation}
The following theorem states the well-posedness of (\ref{weak formul}).
\begin{theorem} \cite{Tobiska_book}
If we assume $c(x,y) -\frac{1}{2}\frac{\partial b_1(x,y)}{\partial x}\geq \gamma_1>0 \;\forall\; (x,y) \in\bar{\Omega}$, then for B(.,.) given in (\ref{weak}) the following holds
    \begin{equation*}
        B(u,u) \geq \frac{1}{2} \|u\|_{1,\varepsilon
        }^2 \quad \forall \; u \in H^1_0(\Omega).
    \end{equation*}
      Henceforth, we denote the positive function $ c - \frac{1}{2} \frac{\partial b_1}{\partial x}$ = $c_0^2$.
      \end{theorem}
Singularly perturbed differential equations are typically characterized by a
small parameter $\varepsilon <<1$ multiplying some or all of the highest order terms in the
differential equation. The solution of such an equation possesses interior or boundary
layers or both. These are very thin regions where the solution and its derivative change
very rapidly. Outside these layers, the solution smoothly approximates that of the lower-order reduced problem (obtained by setting $\varepsilon=0 $ in \eqref{current problem}). In some situations, such as nonlinearities or discontinuous source terms or a sign change in the convection coefficient, interior layers can also be present within the domain. As $\varepsilon\to0$, the solution's derivative within these layers often scales inversely with powers of $\varepsilon$, indicating a significant loss of uniform regularity across the domain. For the problem \eqref{current problem}, turning points in SPPs can be identified by a sign change in the convective coefficient within the domain. At a turning point, the convective term locally vanishes or becomes very small. This means that diffusion, even if small (due to $\varepsilon$), can locally become relatively more dominant, influencing the solution behavior. From both theoretical and computational perspectives, turning point SPPs pose a greater challenge than standard SPPs. Standard numerical methods designed for convection-dominated problems may struggle to accurately resolve both the boundary and interior layers, often leading to spurious oscillations or significant errors if not specifically adapted (e.g., using fitted meshes).\\
To proceed with our analysis, precise bounds on specific derivatives of the solution $u(x,y)$ to the problem \eqref{current problem} are indispensable. This is discussed in the following theorem.
\begin{theorem} \label{decompositionu}
{\cite{Mbayi22},\cite{Tobiska_book}} Let $f \in C^{3,\beta}(\bar{\Omega})$ where $\beta \in [2/3,1]$ and $f$ satisfy the compatibility condition $f(-1, -1) = f(1, 1) = f(1, -1) = f(-1, 1) = 0$, then the exact solution of given problem can be broken down as $u = R + F_1 + F_2 + F_{12}$, where
\begin{align*}
     R =& \text{Regular part of the solution},\\
     F_1 =& \text{Part of the solution containing interior layer near $x=0$},\\
     F_2 =& \text{Part of the solution containing boundary layer near $y=1$ and $y=-1$},\\
     F_{12} = &\text{Part of the solution containing both interior and boundary layers near}\\
     &  \text{the points (0,1) and (0,-1)},
\end{align*}
where for all $(x, y) \in \bar{\Omega}$ and $0 \leq i + j \leq 2$, we have
\begin{align*} 
\left| \frac{\partial^{i+j} R(x,y)}{\partial x^i \partial y^j} \right| &\leq C, \\
\left| \frac{\partial^{i+j} F_1(x,y)}{\partial x^i \partial y^j} \right| &\leq C \varepsilon^{-i} e^{-\frac{\alpha |x| }{\varepsilon}}, \\
\left| \frac{\partial^{i+j} F_2(x,y)}{\partial x^i \partial y^j} \right| &\leq C \varepsilon^{-j/2} \left( e^{-\frac{\beta}{\sqrt{\varepsilon}}(1+y)} + e^{-\frac{\beta}{\sqrt{\varepsilon}}(1-y)} \right), \\
\left| \frac{\partial^{i+j} F_{12}(x,y)}{\partial x^i \partial y^j} \right| &\leq C \varepsilon^{-i-j/2} e^{-\frac{\alpha |x|}{\varepsilon}} \left( e^{-\frac{\beta}{\sqrt{\varepsilon}}(1+y)} + e^{-\frac{\beta}{\sqrt{\varepsilon}}(1-y)} \right).
\end{align*}
\end{theorem}
\section{A classical discretization on a Shishkin mesh}\label{section2}
We denote the $L^2(\Omega)$ norm by $\|.\|$ and the energy norm by $\|.\|_{1,\varepsilon}$ and is defined by
\begin{equation}
        \|v\|^2_{1,\varepsilon} = \varepsilon \|\nabla v\|^2 + \|v\|^2 \quad \forall \; v \in H^1(\Omega)
    \end{equation}
    where $H^1(\Omega)$ is the usual Sobolev space of functions whose first-order weak derivatives lie in $L^2(\Omega)$. We discretize \eqref{weak formul} by means of the standard finite element method on a special rectangular mesh. 
           Let $\mathcal{T}_h$=${\{\tau_i^h\}}_{i=1}^N$ be a family of edge-to-edge rectangle's whose union is $\bar{\Omega}$. Define
        \begin{equation*}
             S_h = \{\mathcal{X} \in \mathcal{C}^0(\bar{\Omega}): \mathcal{X} = 0\; \text{on}\; \partial \Omega, \mathcal{X}|_{\tau_i^h} \; \text{linear in $x$ and $y$}\},
        \end{equation*}
        where $\mathcal{C}^0(\bar{\Omega})$ is space of continuous functions up-to boundary of $\Omega$.
    The standard finite element discretization of \eqref{weak formul} is to find $u_{h}\in S_h$ such that 
         \begin{equation} \label{weak formulation}
             B(u_{h},v_h) = (f,v_h)  \quad \forall \; v_h \in S_h.
         \end{equation}
To effectively address the layering phenomenon inherent in the problem's solution, we employ a layer-adapted mesh of Shishkin type, constructed independently along both the $x$ and $y$ spatial directions. For an even positive integer $N \ge 4$,  the crucial mesh transition parameters in the $x$ and $y$ directions are denoted by $\lambda_x$ and $\lambda_y$, respectively. Note that these newly defined transition parameters are one of the novel contributions to this article. They are defined as:
\[
\lambda_x = min\left(\frac{2\varepsilon}{\alpha} \log\left(\frac{1}{\varepsilon}\right),\frac{1}{2}\right) , \quad \lambda_y = min\left(2\sqrt{\frac{\varepsilon}{\beta}} \log\left(\frac{1}{\varepsilon^{3/2}}\right),  \frac{1}{4}\right).
\]
 Mesh points along the $x$-axis and $y$-axis are,
\begin{equation} \label{eq:3.1}
x_i = 
\begin{cases}
\frac{2\lambda_x}{N} i, & i = 0, 1, \dots, N/2, \\
\lambda_x + \frac{2}{N}(1-\lambda_x)(i - \frac{N}{2}), & i = N/2 + 1, \dots, N,
\end{cases}
\end{equation}
and  \begin{equation} 
y_j = 
\begin{cases}
-1 + \frac{4j}{N}\lambda_y, & j = 0, 1, \dots, \frac{N}{4}, \\
-1 + \lambda_y + \frac{2}{N}(2-2\lambda_y)(j - \frac{N}{4}), & j = \frac{N}{4} + 1, \dots, \frac{3N}{4}, \\
1 - 4\lambda_y(1 - \frac{j}{N}), & j = \frac{3N}{4} + 1, \dots, N.
\end{cases}
\end{equation}
\begin{figure}[ht]
    \centering
    \begin{minipage}{0.44\textwidth}
         \includegraphics[width=1\linewidth]{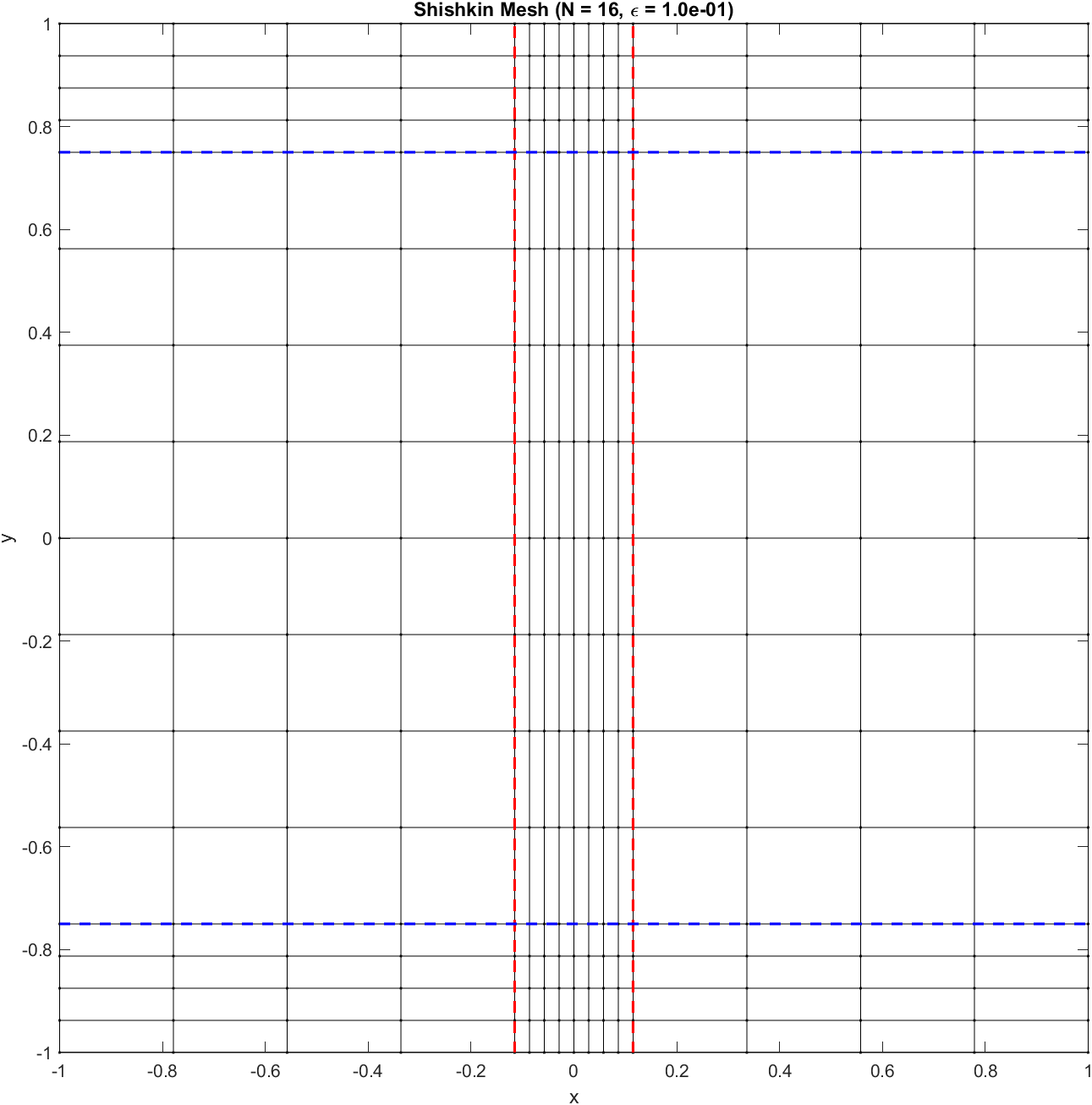}
        \caption{Shishkin mesh}
        \label{fig:img1}
    \end{minipage}
    \hfill
    \begin{minipage}{0.42\textwidth}
        \includegraphics[width=1\linewidth]{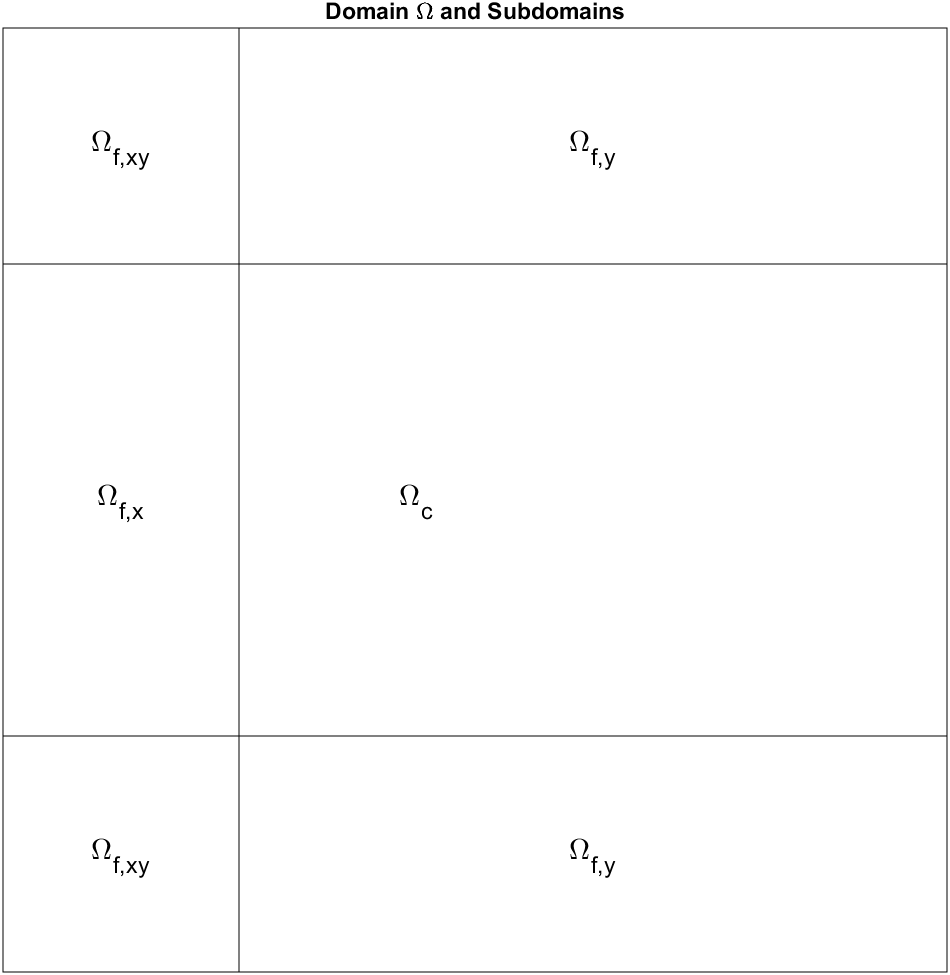} 
        \caption{Domain $\Omega$ and subdomains}
        \label{fig:img2}
    \end{minipage}
\end{figure}
In the $x$-direction, the interval $[0, 1]$ is bisected into two sub-intervals $[0, \lambda_x]$ and $[\lambda_x, 1]$. Each of these sub-intervals is further partitioned such that its closure contains $N/2 + 1$ points, denoted as $X_K^N$. A symmetric partitioning is implicitly applied to $[-1,0]$ to cover the full $x$ range. For the $y$-direction, the interval $[-1, 1]$ is divided into three sub-intervals: $[-1, -1 + \lambda_y]$, $[-1 + \lambda_y, 1 - \lambda_y]$ and $[1 - \lambda_y, 1]$. These sub-intervals are uniformly discretized.  Specifically, their closures are partitioned into $N/4 + 1$, $N/2 + 1$, and $N/4 + 1$ points, respectively, forming the set $Y_J^N$.
The collection of global mesh points is denoted by $\Omega^N$ and defined as the tensor product of these one-dimensional partitions:
\[
\Omega^N = \{(x_i, y_j) \in \Omega : x_i \in X_K^N, y_j \in Y_J^N, i = 0, 1, \dots, N, j = 0, 1, \dots, N\}, \quad N \in \mathbb{N}.
\]
Axiparallel lines passing through the mesh points generate a collection of rectangular mesh elements, denoted as $\mathcal{J}_N$. This subdivision implicitly creates different regions within $\Omega$, often categorized as coarse ($\Omega_c $) and fine ($\Omega_{f,x} $,  $\Omega_{f,y}$, $\Omega_{f,xy}$) mesh region. Note that in the coarse region, the solution is expected to behave nicely, and in the layer regions, the solution is expected to exhibit layers. The mesh widths in the layer regions and coarse region are calculated as:
\begin{align*}
h_x &= 2\lambda_x N^{-1}, \quad \quad H_x = 2(1-\lambda_x)N^{-1},  \\
h_y &= 4\lambda_y N^{-1}, \quad \quad H_y = 4(1-\lambda_y)N^{-1}. 
\end{align*}
Here, $h_x$, $h_y$ denote the mesh widths within the layer regions in the $x$ and $y$ directions, respectively,  while $H_x$, $H_y$ represent the mesh widths outside these regions. Due to the symmetry of the problem with respect to the $y$-axis, we analyze our problem in the reduced domain $\Omega=[0,1]\times[-1,1]$. Within this domain, we perform a decomposition into four non-overlapping sub-regions, tailored to distinguish between areas of expected smooth behavior and those containing layers. This partition is defined as:
         \begin{center}
              $\Omega_c = [\lambda_x,1]\times[-1+\lambda_y,1-\lambda_y], \quad\Omega_{f,x} = [0,\lambda_x]\times[-1+\lambda_y,1-\lambda_y]$,\\
        $\Omega_{f,y} = ([-1,-1+\lambda_y]\cup[1-\lambda_y,1])\times[\lambda_x,1],\quad\Omega_{f,xy} = [0,\lambda_x]\times \left([-1,-1+\lambda_y]\cup[1-\lambda_y,1]\right)$.\\
         \end{center}
         This discretization is followed from \cite{Rajeev_Ranjan}.
\section{Estimates for discrete Green’s function and some auxiliary estimates}\label{section3}
This section is devoted to finding out some important estimates for the discrete Green's function corresponding to the differential operator in \eqref{current problem} and some key auxiliary estimates. These estimates play a fundamental role in deriving the pointwise error estimates. We derive the estimates for the discrete Green's function in $L^2$-norm and energy norm. We note that here the discrete Green's function contains information about the source point. Therefore, these estimates are essential for identifying regions where an $\varepsilon$-independent error estimate is feasible. We shall see that if the source point $(x_m,y_n)\in \Omega_c \cup \Omega_{f,x}$, then we obtain an $\varepsilon$-independent estimate for the discrete Green's function in energy norm and $L^2$-norm, and we can expect to obtain a parameter-uniform $L^\infty$-norm error estimate. In this article, we obtain convergent error estimates specifically for the $x$-layer region, denoted as $\Omega_{f,x}$. However, for the coercive region, we can only establish an upper bound for the error that is independent of $\varepsilon$. This observation is consistent with the findings reported in \cite{Stynes} regarding the coercive region, further validating our analyses.
\begin{lemma}\label{Greens in c} When the source point $(x_m,y_n)  \in\Omega_c$, then $\|G_{mn}\|_{1,\varepsilon}\leq CN.$
\end{lemma}
\textbf{Proof.} From the coercivity of the bilinear form, we have
 \begin{align*}
     C_1 \|G_{mn}\|^2_{1,\varepsilon} &\leq B(G_{mn},G_{mn}) = G_{mn}(x_m,y_n) \leq |G_{mn}(x_m,y_n)|\\
   &  \leq CN^2 \|G_{mn}\|_{L^1(R_{mn})}, \hspace{2cm}  (\text{by inverse estimate})
 \end{align*}
where $R_{mn}$ is the union of those rectangles whose common vertex is $(x_m,y_n)$. Therefore 
\begin{align*}
    C_1 \|G_{mn}\|^2_{1,\varepsilon} &\leq CN^2 \|G_{mn}\|_{L^1(R_{mn})}\\
    &\leq CN^2 |R_{mn}|^{1/2} \|G_{mn}\|_{L^2({R_{mn}})}\\
     &\leq CN^2 \frac{1}{N} \|G_{mn}\|\leq C N \|G_{mn}\|_{{1,\varepsilon}}\\
        \implies\; \|G_{mn}\|_{{1,\varepsilon}} &\leq CN.
\end{align*}
\begin{corollary}
     When $(x_m,y_n) \in \Omega_c$, then $   \|G_{mn}\|\leq C N.$
\end{corollary}
      \begin{lemma}\label{Greens in x} When $(x_m,y_n) \in \Omega_{f,x}$, then $   \|G_{mn}\|_{1,\varepsilon}\leq CN^{1/2} \log^{1/2}(1/\varepsilon).$
      \end{lemma}
    \textbf{Proof.} From the coercivity of the bilinear form, we have 
      \begin{align*}
          C_1\|G_{mn}\|^2_{1,\varepsilon} \leq B(G_{mn},G_{mn})
         & = G_{mn}(x_m,y_n)
         = -\int_{t=x_m}^1G_{mn,x}(t,y_n)\; dt\\
         & \leq \int_{x_m}^1|G_{mn,x}(t,y_n)|\; dt\\
          &\leq C N \int_{y_{n-1}}^{y_n} \int_{x_m}^1|G_{mn,x}(t,s)|\; ds dt\\
          & \leq CN \frac{\varepsilon^{1/2} \log^{1/2}(1/\varepsilon) }{N^{1/2}} \|G_{mn,x}\|\\
          & \leq C N^{1/2} \log^{1/2}(1/\varepsilon) \|G_{mn}\|_{1,\varepsilon}\\
          \implies\;\|G_{mn}\|_{1,\varepsilon}&\leq C  N^{1/2} \log^{1/2} (1/\varepsilon).
      \end{align*}
      \begin{corollary} When $(x_m,y_n) \in \Omega_{f,x}$, then $   \|G_{mn}\|\leq CN^{1/2}\log^{1/2}(1/\varepsilon).$
       \end{corollary} 
\begin{Remark} Note that when $(x_m,y_n)\in\Omega_{f,y}$ we are able to obtain $ \|G_{mn}\|\leq C\varepsilon^{-1/4}\log^{-1/2}(1/\varepsilon)$ and when $(x_m,y_n)\in\Omega_{f,xy}$, then $   \|G_{mn}\|\leq C\varepsilon^{-1/4}$. Due to the $\varepsilon$-dependent estimate for $\|G_{mn}\|$ when $(x_m,y_n)\in \Omega_{f,y}\cup\Omega_{f,xy}$ it is not possible to obtain $\varepsilon$-independent error estimate.  Therefore, it is expected that we do not obtain parameter-independent error estimation when the source point $(x_m,y_n)\in\Omega_{f,y}\cup\Omega_{f,xy}$. This is reflected in the numerical example.
\end{Remark}
\begin{lemma} \label{A_1} If $(x,y) \in \Omega_c\cup\Omega_{f,y}$, then following estimates are holds
\begin{center}
    $\|e^{-\frac{\alpha x}{\varepsilon}}\|_{L^\infty(\Omega_{f,y})} \leq C \varepsilon^\alpha, \quad\quad\quad\; \|e^{-\frac{\alpha x}{\varepsilon}}\|_{L^\infty(\Omega_{c})} \leq C \varepsilon^\alpha$.
\end{center}
\end{lemma} 
\textbf{Proof.} If $(x,y) \in \Omega_{f,y}$, then $x\in[\lambda_x,1]$, $y\in [-1,-1+\lambda_y]\cup [1-\lambda_y,1]$ and when $(x,y) \in \Omega_{c}$, then $x\in[\lambda_x,1]$, $y\in [-1+\lambda_y,1-\lambda_y]$. \\
In the both cases when $x\in [\lambda_x,1]$, we obtain
\begin{align}
    x  \geq \lambda_x & \implies \frac{-\alpha x}{\varepsilon} \leq \frac{-\alpha \lambda_x}{\varepsilon} \nonumber\\
   & \implies e^{\frac{-\alpha x}{\varepsilon}} \leq e^{\frac{-\alpha \lambda_x}{\varepsilon}} = e^{-\frac{\alpha}{\varepsilon}\varepsilon\log\left(\frac{1}{\varepsilon}\right)} 
    = e^{-\alpha \log\left(\frac{1}{\varepsilon}\right)} 
    = e^{ \log\left(\frac{1}{\varepsilon^{-\alpha}}\right)}=\varepsilon^\alpha.\nonumber
\end{align}
This completes the proof.
\begin{corollary} \label{A_3} If $(x,y) \in \Omega_c\cup\Omega_{f,x}$, then following estimates are holds
     $$\|e^{-\frac{\alpha x}{\varepsilon}}\|_{L^\infty(\Omega_{f,x})}=1 \quad \text{and}\quad\|e^{-\frac{\alpha x}{\varepsilon}}\|_{L^\infty(\Omega_{f,xy})}=1.$$
\end{corollary}
\begin{lemma} \label{B_1} If $(x,y) \in \Omega_c\cup\Omega_{f,x}$, then following estimates are holds
\begin{center}
$\|e^{-\frac{\beta}{\sqrt{\varepsilon}}(1+y)}+e^{-\frac{\beta}{\sqrt{\varepsilon}}(1-y)}\|_{L^\infty(\Omega_{c})}\leq C\varepsilon^{\frac{3\beta}{2}},$ \\
 $\|e^{-\frac{\beta}{\sqrt{\varepsilon}}(1+y)}+e^{-\frac{\beta}{\sqrt{\varepsilon}}(1-y)}\|_{L^\infty(\Omega_{f,x})}\leq C \varepsilon^{\frac{3\beta}{2}}.$
\end{center}
\end{lemma} 
\textbf{Proof.} If $(x,y) \in \Omega_{c}$, then $x\in[\lambda_x,1]$, $y\in [-1+\lambda_y,1-\lambda_y]$, we have
\begin{align*}
    \|e^{-\frac{\beta}{\sqrt{\varepsilon}}(1+y)}+e^{-\frac{\beta}{\sqrt{\varepsilon}}(1-y)}\|_{L^\infty(\Omega_{c})} &= \|e^{-\frac{\beta}{\sqrt{\varepsilon}}(1+y)}+e^{-\frac{\beta}{\sqrt{\varepsilon}}(1-y)}\|_{L^\infty(-1+\lambda_y,1-\lambda_y)}\\
    & \leq \|e^{-\frac{\beta}{\sqrt{\varepsilon}}(1+y)}\|_{L^\infty(-1+\lambda_y,1-\lambda_y)} + \|e^{-\frac{\beta}{\sqrt{\varepsilon}}(1-y)}\|_{L^\infty(-1+\lambda_y,1-\lambda_y)},
\end{align*}
and when $(x,y) \in \Omega_{f,x}$, then $x\in[0,\lambda_x]$ and $y\in [-1+\lambda_y,1-\lambda_y]$. We have
\begin{align*}
    \|e^{-\frac{\beta}{\sqrt{\varepsilon}}(1+y)}+e^{-\frac{\beta}{\sqrt{\varepsilon}}(1-y)}\|_{L^\infty(\Omega_{f,x})} &= \|e^{-\frac{\beta}{\sqrt{\varepsilon}}(1+y)}+e^{-\frac{\beta}{\sqrt{\varepsilon}}(1-y)}\|_{L^\infty(-1+\lambda_y,1-\lambda_y)}\\
    & \leq \|e^{-\frac{\beta}{\sqrt{\varepsilon}}(1+y)}\|_{L^\infty(-1+\lambda_y,1-\lambda_y)} + \|e^{-\frac{\beta}{\sqrt{\varepsilon}}(1-y)}\|_{L^\infty(-1+\lambda_y,1-\lambda_y)}.
\end{align*}
For the both cases we need to calculate $\|e^{-\frac{\beta}{\sqrt{\varepsilon}}(1+y)}\|_{L^\infty(-1+\lambda_y,1-\lambda_y)}$. As
\begin{align*}
    y\in[-1+\lambda_y,1-\lambda_y]\implies -1+\lambda_y  \leq y \leq 1-\lambda_y
\end{align*}
First we take $1+y\geq\lambda_y,$ then
\begin{align*}
&-\frac{\beta}{\sqrt{\varepsilon}}(1+y)\leq -\frac{\beta}{\sqrt{\varepsilon}}\lambda_y\implies e^{-\frac{\beta}{\sqrt{\varepsilon}}(1+y)}\leq e^{-\frac{\beta}{\sqrt{\varepsilon}}\lambda_y}\\
  & \implies e^{-\frac{\beta}{\sqrt{\varepsilon}}(1+y)}\leq e^{-\frac{\beta}{\sqrt{\varepsilon}}\sqrt{\varepsilon}\log\left(\frac{1}{\varepsilon^{3/2}}\right)} = e^{-\beta \log\left(\frac{1}{\varepsilon^{3/2}}\right)} 
     = e^{\log\left(\frac{1}{\varepsilon^{\frac{-3\beta}{2}}}\right)} 
    = \varepsilon^{\frac{3\beta}{2}},
\end{align*}
and if $ y\leq 1-\lambda_y,$     then
\begin{align*}
 &-\frac{\beta}{\sqrt{\varepsilon}}(1-y)\leq -\frac{\beta}{\sqrt{\varepsilon}}\lambda_y
\implies e^{-\frac{\beta}{\sqrt{\varepsilon}}(1-y)}\leq e^{-\frac{\beta}{\sqrt{\varepsilon}}\lambda_y}\\
   &\implies e^{-\frac{\beta}{\sqrt{\varepsilon}}(1-y)}\leq e^{-\frac{\beta}{\sqrt{\varepsilon}}\sqrt{\varepsilon}\log\left(\frac{1}{\varepsilon^{3/2}}\right)}  = e^{-\beta \log\left(\frac{1}{\varepsilon^{3/2}}\right)} 
     = e^{\log\left(\frac{1}{\varepsilon^{\frac{-3\beta}{2}}}\right)} 
    = \varepsilon^{\frac{3\beta}{2}}.
\end{align*}
This completes the proof.
\begin{corollary} \label{B_3} If $(x,y) \in \Omega_{f,y}\cup\Omega_{f,xy}$, then following estimates are holds
\begin{center}
    $\|e^{-\frac{\beta}{\sqrt{\varepsilon}}(1+y)}+e^{-\frac{\beta}{\sqrt{\varepsilon}}(1-y)}\|_{L^\infty(\Omega_{f,y})}=1, $
and $\|e^{-\frac{\beta}{\sqrt{\varepsilon}}(1+y)}+e^{-\frac{\beta}{\sqrt{\varepsilon}}(1-y)}\|_{L^\infty(\Omega_{f,xy})}=1. $
\end{center}
\end{corollary}
\section{Estimates related to diffusion terms}\label{section4}
In this section, we aim to bound the diffusion portion of the bilinear form B(.,.) with its arguments parts of interpolation error $(u-Iu)$ and the discrete Green's function $ G_{mn}.$ We start by presenting the following auxiliary results.
\begin{theorem}\label{ajeev}
\cite{Rajeev_Ranjan} Assume that $K \in \mathcal{J}_N$ is a mesh rectangle with sides that are perpendicular to the coordinate axes. Assume that $u \in H^3(K)$ and that $Iu$ is its nodal bilinear interpolate on $K$. Hence, we obtain for any bilinear function $v^N$ defined on $K$, then 
\begin{align*}
    \left| \int_K (Iu - u)_x v_x^N \,dxdy \right| &\le C h_{y,K}^2 \|u_{xyy}\|_{L^2(K)} \|v_x^N\|_{L^2(K)}, \\
    \left| \int_K (Iu - u)_y v_y^N \,dxdy \right| &\le C h_{x,K}^2 \|u_{xxy}\|_{L^2(K)} \|v_y^N\|_{L^2(K)}.
\end{align*}
\end{theorem}
\begin{lemma} \label{diffusion R} If source point $(x_m,y_n) \in\Omega_c\cup \Omega_{f,x}$, then $|\varepsilon (\nabla(R-IR), \nabla G_{mn})|  \leq \frac{C(\log(1/\varepsilon^{3/2}))^{3}}{N}.$
     \end{lemma} 
\textbf{Proof.} First we take the source point $(x_m,y_n)\in \Omega_{f,x}$ and decompose the domain in four parts, we have
\begin{align}
    |\varepsilon (\nabla(R-IR),\nabla G_{mn})|
    & \leq |\varepsilon (\nabla(R-IR),\nabla G_{mn})_{\Omega_{f,xy}}| 
     + |\varepsilon (\nabla(R-IR),\nabla G_{mn})_{\Omega_{f,x}}|  \nonumber\\
    &+ |\varepsilon (\nabla(R-IR),\nabla G_{mn})_{\Omega_{f,y}}|  + |\varepsilon (\nabla(R-IR),\nabla G_{mn})_{\Omega_{c}}|.   \label{eq:R}
\end{align}
Now, estimate the first term of the R.H.S. of  (\ref{eq:R}) in the $x$-direction   is 
  $|\varepsilon ((R-IR)_x,G_{mn,x})_{\Omega_{f,xy}}|$. By using Theorem \ref{decompositionu}, Theorem \ref{ajeev} and furthermore, we obtain
\begin{align}
    |\varepsilon ((R-IR)_x, G_{mn,x})_{\Omega_{f,xy}}| 
    & \leq  C \sqrt{\varepsilon} \; h_y^2 \; \|{R}_{xyy}\|_{\Omega_{f,xy}}  \; \sqrt{\varepsilon}\|G_{mn,x}\|_{{\Omega_{f,xy}} } \nonumber  \\
    & \leq C \sqrt{\varepsilon} \left(\frac{\sqrt{\varepsilon}\log(1/\varepsilon^{3/2})}{N} \right)^2 \sqrt{|\Omega_{f,xy}|} \:\|G_{mn}\|_{1,\varepsilon} \nonumber\\
    &\leq C \sqrt{\varepsilon} \left(\frac{\sqrt{\varepsilon}\log(1/\varepsilon^{3/2})}{N} \right)^2\sqrt{\varepsilon \log(1/\varepsilon)} \sqrt{\sqrt{\varepsilon}\log(1/\varepsilon^{3/2})}\; \sqrt{N \log\left(\frac{1}{\varepsilon}\right)}\nonumber\\
    &\leq C \frac{\varepsilon^{9/4}(\log(1/\varepsilon^{3/2}))^{7/2}}{N^{3/2}} \leq C \frac{(\log(1/\varepsilon^{3/2}))^{7/2}}{N^{3/2}}.
\end{align}
 Now, estimate the first term of the R.H.S. of  (\ref{eq:R}) in the $y$-direction   is
  $|\varepsilon ((R-IR)_y,G_{mn,y})_{\Omega_{f,xy}}|$. By using Theorem \ref{decompositionu}, Theorem \ref{ajeev} and furthermore, we obtain
\begin{align}
    |\varepsilon ((R-IR)_y, G_{mn,y})_{\Omega_{f,xy}}| 
    & \leq  C \sqrt{\varepsilon} \; h_x^2 \; \|{F_1}_{yxx}\|_{\Omega_{f,xy}} \; \sqrt{\varepsilon}\|G_{mn,y}\|_{{\Omega_{f,xy}} } \nonumber  \\
    & \leq C \sqrt{\varepsilon} \left(\frac{\varepsilon\log(1/\varepsilon)}{N} \right)^2 \sqrt{|\Omega_{f,xy}|} \;\|G_{mn}\|_{1,\varepsilon}\nonumber\\
     &\leq C \sqrt{\varepsilon}  \left(\frac{\varepsilon\log(1/\varepsilon)}{N} \right)^2 \sqrt{\varepsilon \log(1/\varepsilon)} \sqrt{\sqrt{\varepsilon}\log(1/\varepsilon^{3/2})} \; \sqrt{N \log\left(\frac{1}{\varepsilon}\right)}\nonumber\\
    &\leq C \frac{\varepsilon^{13/4}(\log(1/\varepsilon^{3/2}))^{7/2}}{N^{3/2}} \leq C \frac{(\log(1/\varepsilon^{3/2}))^{7/2}}{N^{3/2}}.
\end{align}
Now, estimate the second term of the R.H.S. of  (\ref{eq:R}) in the $x$-direction   is
  $|\varepsilon ((R-IR)_x,G_{mn,x})_{\Omega_{f,x}}|$. By using \ref{decompositionu}, Theorem \ref{ajeev} and furthermore, we obtain
\begin{align}
    |\varepsilon ((R-IR)_x, G_{mn,x})_{\Omega_{f,x}}| 
    & \leq  C \sqrt{\varepsilon} \; h_y^2 \; \|{R}_{xyy}\|_{\Omega_{f,x}} \; \sqrt{\varepsilon}\|G_{mn,x}\|_{{\Omega_{f,x}} } \nonumber  \\
    & \leq C \sqrt{\varepsilon} \left(\frac{1}{N} \right)^2 \sqrt{|\Omega_{f,x}|} \:\|G_{mn}\|_{1,\varepsilon}\nonumber\\
    &\leq C \frac{\varepsilon\log(1/\varepsilon)}{N^{3/2}} \leq C \frac{\log(1/\varepsilon)^{1/2}}{N^{3/2}}.
\end{align}
 Now, estimate the second term of the R.H.S. of  (\ref{eq:R}) in the $y$-direction   is
  $|\varepsilon ((R-IR)_y,G_{mn,y})_{\Omega_{f,x}}|$. By using Theorem \ref{decompositionu}, Theorem \ref{ajeev} and furthermore, we obtain
\begin{align}
    |\varepsilon ((R-IR)_y, G_{mn,y})_{\Omega_{f,x}}| 
    & \leq  C \sqrt{\varepsilon} \; h_x^2 \; \|{F_1}_{yxx}\|_{\Omega_{f,x}} \; \sqrt{\varepsilon}\|G_{mn,y}\|_{{\Omega_{f,x}} } \nonumber  \\
    & \leq C \sqrt{\varepsilon} \left(\frac{\varepsilon\log(1/\varepsilon)}{N} \right)^2 \sqrt{|\Omega_{f,x}|} \;\|G_{mn}\|_{1,\varepsilon}\nonumber\\
     &\leq C \sqrt{\varepsilon}  \left(\frac{\varepsilon\log(1/\varepsilon)}{N} \right)^2 \sqrt{\varepsilon \log(1/\varepsilon)} \; \sqrt{N \log\left(\frac{1}{\varepsilon}\right)}\nonumber\\
    &\leq C \frac{\varepsilon^3(\log(1/\varepsilon))^3}{N^{3/2}} \leq C \frac{(\log(1/\varepsilon))^3}{N^{3/2}}.
\end{align}
Now, estimate the third term of the R.H.S. of  (\ref{eq:R}) in the $x$-direction   is
  $|\varepsilon ((R-IR)_x,G_{mn,x})_{\Omega_{f,y}}|$. By using Theorem \ref{decompositionu}, Theorem \ref{ajeev} and furthermore, we obtain
\begin{align}
    |\varepsilon ((R-IR)_x, G_{mn,x})_{\Omega_{f,y}}| 
    & \leq  C \sqrt{\varepsilon} \; h_y^2 \; \|{R}_{xyy}\|_{\Omega_{f,y}} \; \sqrt{\varepsilon}\|G_{mn,x}\|_{{\Omega_{f,y}} } \nonumber  \\
    & \leq C \sqrt{\varepsilon} \left(\frac{\sqrt{\varepsilon}\log(1/\varepsilon^{3/2})}{N} \right)^2 \sqrt{|\Omega_{f,y}|} \:\|G_{mn}\|_{1,\varepsilon}\nonumber\\
    &\leq C \sqrt{\varepsilon} \left(\frac{\sqrt{\varepsilon}\log(1/\varepsilon^{3/2})}{N} \right)^2\sqrt{\sqrt{\varepsilon}\log(1/\varepsilon^{3/2})} \; \sqrt{N \log\left(\frac{1}{\varepsilon}\right)}\nonumber\\
    &\leq C \frac{\varepsilon^{7/4}(\log(1/\varepsilon^{3/2}))^3}{N^{3/2}} \leq C \frac{(\log(1/\varepsilon^{3/2}))^3}{N^{3/2}}.
\end{align}
 Now, estimate the third term of the R.H.S. of  (\ref{eq:R}) in the $y$-direction   is
  $|\varepsilon ((R-IR)_y,G_{mn,y})_{\Omega_{f,y}}|$. By using Theorem \ref{decompositionu}, Theorem \ref{ajeev} and furthermore, we obtain
\begin{align}
    |\varepsilon ((R-IR)_y, G_{mn,y})_{\Omega_{f,y}}| 
    & \leq  C \sqrt{\varepsilon} \; h_x^2 \; \|{F_1}_{yxx}\|_{\Omega_{f,y}} \; \sqrt{\varepsilon}\|G_{mn,y}\|_{{\Omega_{f,y}} } \nonumber  \\
    & \leq C \sqrt{\varepsilon} \left(\frac{1}{N} \right)^2 \sqrt{|\Omega_{f,y}|} \;\|G_{mn}\|_{1,\varepsilon}\nonumber\\
     &\leq C \sqrt{\varepsilon}  \left(\frac{1}{N} \right)^2 \sqrt{\sqrt{\varepsilon}\log(1/\varepsilon^{3/2})} \; \sqrt{N \log\left(\frac{1}{\varepsilon}\right)}\nonumber\\
    &\leq C \frac{\varepsilon^{3/4}\log(1/\varepsilon^{3/2})}{N^{3/2}} \leq C \frac{\log(1/\varepsilon^{3/2})}{N^{3/2}}.
\end{align}
 Now, estimate the fourth term of the R.H.S. of  (\ref{eq:R}) 
  $|\varepsilon (\nabla(R-IR),\nabla G_{mn})_{\Omega_{c}}| $. By using Theorem \ref{decompositionu}, Theorem \ref{ajeev} and furthermore, we obtain
\begin{align}
    |\varepsilon (\nabla(R-IR), \nabla G_{mn})_{\Omega_{c}}|  &\leq \varepsilon^{1/2}\|\nabla(R-IR)\| \varepsilon^{1/2}\|G_{mn,x}\|_{\Omega_c} \nonumber\\
    & \leq   \varepsilon^{1/2} CN^{-2} \sqrt{N \log\left(\frac{1}{\varepsilon}\right)}  \nonumber  \\
    &\leq C \varepsilon^{1/2}\frac{(\log(1/\varepsilon))^{1/2}}{N^{3/2}} \leq \frac{C (\log(1/\varepsilon))^{1/2}}{N^{3/2}}. \label{eq:11}
\end{align}
Applying similar arguments we can show that when source point $(x_m,y_n) \in \Omega_{c}$, then we have
\begin{center}
    $|\varepsilon (\nabla(R-IR),\nabla G_{mn})| \leq C \frac{(\log(1/\varepsilon^{3/2})^3}{N}.$
\end{center}
Hence, for the source point $(x_m,y_n) \in \Omega_{c}\cup \Omega_{f,x}$, we have 
$$|\varepsilon (\nabla(R-IR), \nabla G_{mn})|\leq \frac{C(\log(1/\varepsilon^{3/2}))^{3}}{N}$$
\begin{lemma} \label{diffusion F_1}  If the source point $(x_m,y_n)\in\Omega_c\cup \Omega_{f,x}$, then for $\alpha\geq 3/2$ following estimate hold 
\begin{center}
$|\varepsilon (\nabla(F_1-IF_1), \nabla G_{mn})| \leq C  \frac{(\log(1/\varepsilon^{3/2}))^{3}}{N}.$    
\end{center}
\end{lemma} 
\textbf{Proof.} First we take the source point $(x_m,y_n)\in \Omega_{f,x}$ and decompose the domain in four parts, we have
\begin{align}
    |\varepsilon (\nabla(F_1-IF_1),\nabla G_{mn})|
    & \leq |\varepsilon (\nabla(F_1-IF_1),\nabla G_{mn})_{\Omega_{f,xy}}|  
     + |\varepsilon (\nabla(F_1-IF_1),\nabla G_{mn})_{\Omega_{f,x}}|  \nonumber\\
    &+ |\varepsilon (\nabla(F_1-IF_1),\nabla G_{mn})_{\Omega_{f,y}}|  + |\varepsilon (\nabla(F_1-IF_1),\nabla G_{mn})_{\Omega_{c}}|   \label{eq:13}
\end{align}
Now, estimate the first term of the R.H.S. of  (\ref{eq:13}) in the $x$-direction   is
  $|\varepsilon ((F_1-IF_1)_x,G_{mn,x})_{\Omega_{f,xy}}|$. By using Theorem \ref{decompositionu}, Theorem \ref{ajeev}, Corollary \ref{A_3} and furthermore, we obtain
\begin{align}
    |\varepsilon ((F_1-IF_1)_x, G_{mn,x})_{\Omega_{f,xy}}| 
    & \leq  C \sqrt{\varepsilon} \; h_y^2 \; \|{F_1}_{xyy}\|_{\Omega_{f,xy}} \; \sqrt{\varepsilon}\|G_{mn,x}\|_{{\Omega_{f,xy}} } \nonumber  \\
    & \leq C \sqrt{\varepsilon} \left(\frac{\sqrt{\varepsilon}\log(1/\varepsilon^{3/2})}{N} \right)^2 \sqrt{|\Omega_{f,xy}|} \frac{1}{\varepsilon} \|e^{-\frac{\alpha x}{\varepsilon}}\|_{L^\infty({\Omega_{f,xy}})} \|G_{mn}\|_{1,\varepsilon}\nonumber\\
    &\leq C\sqrt{\varepsilon} \left(\frac{\log(1/\varepsilon^{3/2})}{N} \right)^2 \sqrt{\varepsilon \log(1/\varepsilon)} \sqrt{\sqrt{\varepsilon}\log(1/\varepsilon^{3/2})} \;\sqrt{N \log\left(\frac{1}{\varepsilon}\right)} \nonumber\\
    &\leq C  \varepsilon^{5/4} \frac{(\log(1/\varepsilon^{3/2}))^{7/2}}{N^{3/2}}  \leq C  \frac{(\log(1/\varepsilon^{3/2}))^{7/2}}{N^{3/2}}.
\end{align}
 Now, estimate the first term of the R.H.S. of  (\ref{eq:13}) in the $y$-direction   is
  $|\varepsilon ((F_1-IF_1)_y,G_{mn,y})_{\Omega_{f,xy}}|$. By using Theorem \ref{decompositionu}, Theorem \ref{ajeev}, Corollary \ref{A_3} and furthermore, we obtain
\begin{align}
    |\varepsilon ((F_1-IF_1)_y, G_{mn,y})_{\Omega_{f,xy}}| 
    & \leq  C \sqrt{\varepsilon} \; h_x^2 \; \|{F_1}_{yxx}\|_{\Omega_{f,xy}} \; \sqrt{\varepsilon}\|G_{mn,y}\|_{{\Omega_{f,xy}} }  \nonumber  \\
    & \leq C \sqrt{\varepsilon} \left(\frac{\varepsilon\log(1/\varepsilon)}{N} \right)^2 \sqrt{|\Omega_{f,xy}|} \frac{1}{\varepsilon^2} \|e^{-\frac{\alpha x}{\varepsilon}}\|_{L^\infty({\Omega_{f,xy}})} \|G_{mn}\|_{1,\varepsilon}\nonumber\\
    &\leq C\sqrt{\varepsilon} \left(\frac{\log(1/\varepsilon^{3/2})}{N} \right)^2 \sqrt{\varepsilon \log(1/\varepsilon)} \sqrt{\sqrt{\varepsilon}\log(1/\varepsilon^{3/2})} \; \sqrt{N \log\left(\frac{1}{\varepsilon}\right)}\nonumber\\
    &\leq C  \varepsilon^{5/4} \frac{(\log(1/\varepsilon^{3/2}))^{7/2}}{N^{3/2}} \leq C  \frac{(\log(1/\varepsilon^{3/2}))^{7/2}}{N^{3/2}}.
\end{align}
 Now, estimate the second term of the R.H.S. of  (\ref{eq:13}) in the $x$-direction
  $|\varepsilon ((F_1-IF_1)_x,G_{mn,x})_{\Omega_{f,x}}|$. By using Theorem \ref{decompositionu}, Theorem \ref{ajeev}, Corollary \ref{A_1} and furthermore, we obtain
\begin{align}
    |\varepsilon ((F_1-IF_1)_x, G_{mn,x})_{\Omega_{f,x}}| 
    & \leq  C \sqrt{\varepsilon} \; h_y^2 \; \|{F_1}_{xyy}\|_{\Omega_{f,x}} \; \sqrt{\varepsilon}\|G_{mn,x}\|_{{\Omega_{f,x}} } \nonumber  \\
    & \leq C \sqrt{\varepsilon} \frac{1}{N^2} \sqrt{|\Omega_{f,x}|} \frac{1}{\varepsilon} \|e^{-\frac{\alpha x}{\varepsilon}}\|_{L^\infty({\Omega_{f,x}})} \|G_{mn}\|_{1,\varepsilon}\nonumber\\
    &\leq C\sqrt{\varepsilon}\frac{1}{\varepsilon N^2}\sqrt{\varepsilon \log(1/\varepsilon)} \; \sqrt{N \log\left(\frac{1}{\varepsilon}\right)} \nonumber\\
    &\leq C   \frac{\log(1/\varepsilon)}{N^{3/2}}.
\end{align}
 Now, estimate the second term of the R.H.S. of  (\ref{eq:13}) in the $y$-direction is
  $|\varepsilon ((F_1-IF_1)_y,G_{mn,y})_{\Omega_{f,x}}|$. By using Theorem \ref{decompositionu}, Theorem \ref{ajeev}, Corollary \ref{A_3} and furthermore, we obtain
\begin{align}
    |\varepsilon ((F_1-IF_1)_y, G_{mn,y})_{\Omega_{f,x}}| 
    & \leq  C \sqrt{\varepsilon} \; h_x^2 \; \|{F_1}_{yxx}\|_{\Omega_{f,x}} \; \sqrt{\varepsilon}\|G_{mn,y}\|_{{\Omega_{f,x}} } \nonumber  \\
    & \leq C \sqrt{\varepsilon} \left(\frac{\varepsilon\log(1/\varepsilon)}{N} \right)^2 \sqrt{|\Omega_{f,x}|} \frac{1}{\varepsilon^2} \|e^{-\frac{\alpha x}{\varepsilon}}\|_{L^\infty({\Omega_{f,x}})} \|G_{mn}\|_{1,\varepsilon}\nonumber\\
    &\leq C\sqrt{\varepsilon} \left(\frac{\log(1/\varepsilon)}{N} \right)^2 \sqrt{\varepsilon \log(1/\varepsilon)}  \; \sqrt{N \log\left(\frac{1}{\varepsilon}\right)} \nonumber\\
    &\leq C  \varepsilon^{1/2} \frac{(\log(1/\varepsilon))^3}{N^{3/2}} \leq C  \frac{(\log(1/\varepsilon))^3}{N^{3/2}}.
\end{align}
 Now, estimate the third term of the R.H.S. of  (\ref{eq:13}) in the $x$-direction is
  $|\varepsilon ((F_1-IF_1)_x,G_{mn,x})_{\Omega_{f,y}}|$. By using Theorem \ref{decompositionu}, Theorem \ref{ajeev}, Lemma \ref{A_1} and furthermore, we obtain
\begin{align}
    |\varepsilon ((F_1-IF_1)_x, G_{mn,x})_{\Omega_{f,y}}| 
    & \leq  C \sqrt{\varepsilon} \; h_y^2 \; \|{F_1}_{xyy}\|_{\Omega_{f,y}} \; \sqrt{\varepsilon}\|G_{mn,x}\|_{{\Omega_{f,y}} } \nonumber  \\
    & \leq C \sqrt{\varepsilon} \left(\frac{\sqrt{\varepsilon}\log(1/\varepsilon^{3/2})}{N} \right)^2 \sqrt{|\Omega_{f,y}|} \frac{1}{\varepsilon} \|e^{-\frac{\alpha x}{\varepsilon}}\|_{L^\infty({\Omega_{f,y}})} \|G_{mn}\|_{1,\varepsilon}\nonumber\\
    &\leq C\sqrt{\varepsilon}\frac{(\log(1/\varepsilon^{3/2}))^2}{ N^2}\sqrt{\sqrt{\varepsilon} \log(1/\varepsilon^{3/2})}\;\varepsilon^{\alpha} \;\sqrt{N \log\left(\frac{1}{\varepsilon}\right)} \nonumber\\
    & \leq C \varepsilon^{\alpha+3/4}  \frac{(\log(1/\varepsilon^{3/2}))^{3}}{N^{3/2}}\leq C  \frac{(\log(1/\varepsilon^{3/2}))^{3}}{N^{3/2}}.
\end{align}
Now, estimate the third term of the R.H.S. of  (\ref{eq:13}) in the $y$-direction is
  $|\varepsilon ((F_1-IF_1)_y,G_{mn,y})_{\Omega_{f,y}}|$. By using Theorem \ref{decompositionu}, Theorem \ref{ajeev}, Lemma \ref{A_1} and furthermore, we obtain
  \begin{align}
    |\varepsilon ((F_1-IF_1)_y, G_{mn,y})_{\Omega_{f,y}}| 
    & \leq  C \sqrt{\varepsilon} \; h_x^2 \; \|{F_1}_{yxx}\|_{\Omega_{f,y}} \; \sqrt{\varepsilon}\|G_{mn,y}\|_{{\Omega_{f,y}} } \nonumber  \\
    & \leq C \sqrt{\varepsilon} \frac{1}{N^2} \sqrt{|\Omega_{f,y}|} \frac{1}{\varepsilon^2} \|e^{-\frac{\alpha x}{\varepsilon}}\|_{L^\infty({\Omega_{f,y}})} \|G_{mn}\|_{1,\varepsilon} \nonumber\\
    &\leq C\sqrt{\varepsilon} \frac{1}{N^2} \sqrt{\sqrt{\varepsilon} \log(1/\varepsilon^{3/2})} \frac{1}{\varepsilon^2} \|e^{-\frac{\alpha x}{\varepsilon}}\|_{L^\infty({\Omega_{f,y}})}\; \sqrt{N \log\left(\frac{1}{\varepsilon}\right)}\nonumber\\
    &\leq C   \frac{\log(1/\varepsilon^{3/2})}{N^{3/2}} \varepsilon^{\alpha-3/2}  \leq C \frac{\log(1/\varepsilon^{3/2})}{N^{3/2}} \quad  \text{when} \; \alpha\geq 3/2.
\end{align}
 Now, estimate the fourth  term of the R.H.S. of (\ref{eq:13}) in the $x$-direction   is 
  $|\varepsilon ((F_1-IF_1)_x,G_{mn,x})_{\Omega_{c}}|$. By using Theorem \ref{decompositionu}, Theorem \ref{ajeev}, Lemma \ref{A_1} and furthermore, we obtain
 \begin{align}
    |\varepsilon ((F_1-IF_1)_x, G_{mn,x})_{\Omega_{c}}|
    & \leq  C \sqrt{\varepsilon} \; h_y^2 \; \|{F_1}_{xyy}\|_{\Omega_{c}} \; \sqrt{\varepsilon}\|G_{mn,x}\|_{{\Omega_{c}} } \nonumber  \\
    & \leq C \sqrt{\varepsilon} \left(\frac{1}{N} \right)^2 \sqrt{|\Omega_{c}|} \frac{1}{\varepsilon} \|e^{-\frac{\alpha x}{\varepsilon}}\|_{L^\infty({\Omega_{c}})} \|G_{mn}\|_{1,\varepsilon}\nonumber\\
    &\leq C \frac{1}{N^2} \frac{1}{\sqrt{\varepsilon}} \|e^{-\frac{\alpha x}{\varepsilon}}\|_{L^\infty({\Omega_{c}})} \sqrt{N \log\left(\frac{1}{\varepsilon}\right)} \nonumber\\
    &\leq C  \varepsilon^{-1/2} \frac{(\log(1/\varepsilon))^{1/2}}{N^{3/2}}\varepsilon^{\alpha} \leq  \frac{C (\log(1/\varepsilon))^{1/2}}{N^{3/2}}.  \quad \text{when} \; \alpha\geq 1/2.
\end{align}
Now, estimate the fourth  term of the R.H.S. of (\ref{eq:13})  in the $y$-direction   is 
  $|\varepsilon ((F_1-IF_1)_y,G_{mn,y})_{\Omega_{c}}|$. By using Theorem \ref{decompositionu}, Theorem \ref{ajeev}, Lemma \ref{A_1} and furthermore, we obtain
 \begin{align}
    |\varepsilon ((F_1-IF_1)_y, G_{mn,y})_{\Omega_{c}}| 
    & \leq  C \sqrt{\varepsilon} \; h_x^2 \; \|{F_1}_{yxx}\|_{\Omega_{c}} \; \sqrt{\varepsilon}\|G_{mn,y}\|_{{\Omega_{c}} } \nonumber  \\
    & \leq C \sqrt{\varepsilon} \left(\frac{1}{N} \right)^2 \sqrt{|\Omega_{c}|} \frac{1}{\varepsilon^2} \|e^{-\frac{\alpha x}{\varepsilon}}\|_{L^\infty({\Omega_{c}})} \|G_{mn}\|_{1,\varepsilon} \nonumber\\
    &\leq C \frac{1}{N^2} \frac{1}{\varepsilon^{3/2}} \|e^{-\frac{\alpha x}{\varepsilon}}\|_{L^\infty({\Omega_{c}})} \sqrt{N \log\left(\frac{1}{\varepsilon}\right)} \nonumber\\
    &  \leq  C \frac{\varepsilon^{\alpha-3/2}}{N^{3/2}} (\log(1/\varepsilon))^{1/2}
    \leq \frac{C(\log(1/\varepsilon))^{1/2}}{N^{3/2}}. \quad  \text{when} \; \alpha\geq 3/2.
\end{align}
Applying similar arguments we can show that when source point $(x_m,y_n) \in \Omega_{c}$, then we have
$$|\varepsilon (\nabla(F_1-IF_1),\nabla G_{mn})| \leq C \frac{(\log(1/\varepsilon^{3/2})^3}{N}.$$
Hence, when source point $(x_m,y_n) \in \Omega_{c}\cup \Omega_{f,x}$, we have 
$$|\varepsilon (\nabla(F_1-IF_1),\nabla G_{mn})| \leq C \frac{(\log(1/\varepsilon^{3/2})^{3}}{N}.$$
\begin{lemma} \label{diffusion F_2}
If the source point $(x_m,y_n)\in\Omega_c\cup \Omega_{f,x}$, then for $\beta\geq 1/3$ following estimate hold  $$|\varepsilon (\nabla(F_2-IF_2), \nabla G_{mn})| \leq C  \frac{(\log(1/\varepsilon^{3/2}))^{3}}{N}.$$
\end{lemma} 
\textbf{Proof.} First we take the source point $(x_m,y_n)\in \Omega_{f,x}$ and decompose the domain in four parts, we have
\begin{align}
   |\varepsilon (\nabla(F_2-IF_2),\nabla G_{mn})|
    & \leq |\varepsilon (\nabla(F_2-IF_2),\nabla G_{mn})_{\Omega_{f,xy}}|  
     + |\varepsilon (\nabla(F_2-IF_2),\nabla G_{mn})_{\Omega_{f,x}}|  \nonumber\\
    &+ |\varepsilon (\nabla(F_2-IF_2),\nabla G_{mn})_{\Omega_{f,y}}|  
    + |\varepsilon (\nabla(F_2-IF_2),\nabla G_{mn})_{\Omega_{c}}|    \label{F_2}
\end{align}
 Now, estimate the first term of the R.H.S. of  (\ref{F_2}) in the $x$-direction   is
  $|\varepsilon ((F_2-IF_2)_x,G_{mn,x})_{\Omega_{f,xy}}|$. By using Theorem \ref{decompositionu}, Theorem \ref{ajeev}, Corollary \ref{B_3} and furthermore, we obtain
\begin{align}
    |\varepsilon ((F_2-IF_2)_x, G_{mn,x})_{\Omega_{f,xy}}| &\leq \sqrt{\varepsilon}\|(F_2-IF_2)_x\|_{\Omega_{f,xy}} \sqrt{\varepsilon}\|G_{mn,x}\|_{\Omega_{f,xy}}  \nonumber\\
    & \leq  C \varepsilon^{1/2} \; h_y^2 \|{F_2}_{xyy}\|_{\Omega_{f,xy}} \; 
 \; \|G_{mn}\|_{{_{1,\varepsilon}} }\nonumber  \\
    &\leq C \varepsilon^{1/2} \left(\frac{\sqrt{\varepsilon}\log(1/\varepsilon^{3/2})}{N} \right)^2 \frac{1}{\varepsilon} \|e^{-\frac{\beta}{\sqrt{\varepsilon}}(1+y)}+e^{-\frac{\beta}{\sqrt{\varepsilon}}(1-y)}\|_{L^\infty(\Omega_{f,xy})} \nonumber\\
   & \leq C \varepsilon^{1/2} \frac{(\log (1/\varepsilon^{3/2})^{2}}{N} \varepsilon^{3/4} \log (1/\varepsilon^{3/2}) \sqrt{N \log\left(\frac{1}{\varepsilon}\right)} \nonumber\\
    &\leq C \varepsilon^{5/4} \frac{(\log (1/\varepsilon^{3/2}))^{7/2}}{N^{3/2}} 
     \leq C \frac{(\log (1/\varepsilon^{3/2}))^{7/2}}{N^{3/2}}.
\end{align}
Now, estimate the first term of the R.H.S. of  (\ref{F_2}) in the $y$- direction is
  $|\varepsilon ((F_2-IF_2)_y,G_{mn,y})_{\Omega_{f,xy}}|$. By using Theorem \ref{decompositionu}, Theorem \ref{ajeev}, Corollary \ref{B_3} and furthermore, we obtain
\begin{align}
    |\varepsilon ((F_2-IF_2)_y, G_{mn,y})_{\Omega_{f,xy}}| &\leq \sqrt{\varepsilon}\|(F_2-IF_2)_y\|_{\Omega_{f,y}} \sqrt{\varepsilon}\|G_{mn,y}\|_{\Omega_{f,xy}} \nonumber\\
    & \leq  C \varepsilon^{1/2} \; h_x^2 \|{F_2}_{yxx}\|_{\Omega_{f,xy}} \; 
 \; \|G_{mn}\|_{{_{1,\varepsilon}} }\nonumber  \\
    &\leq  \frac{C \varepsilon^{1/2}}{N^2\varepsilon^{1/2}} \|e^{-\frac{\beta}{\sqrt{\varepsilon}}(1+y)}+e^{-\frac{\beta}{\sqrt{\varepsilon}}(1-y)}\|_{L^\infty(\Omega_{f,y})}  \sqrt{|\Omega_{f,xy}|}\;\sqrt{N \log\left(\frac{1}{\varepsilon}\right)} \nonumber\\
    &\leq C\varepsilon^{3/4} \frac{(\log (1/\varepsilon^{3/2}))^{3/2}}{N^{3/2}} 
    \leq C  \frac{(\log (1/\varepsilon^{3/2}))^{3/2}}{N^{3/2}}.
\end{align}
 Now, estimate second term of the R.H.S. of  (\ref{F_2}) in the $x$-direction is
  $|\varepsilon ((F_2-IF_2)_x,G_{mn,x})_{\Omega_{f,x}}|$. By using Theorem \ref{decompositionu}, Theorem \ref{ajeev}, Lemma \ref{B_1} and furthermore, we obtain
\begin{align}
    |\varepsilon ((F_2-IF_2)_x, G_{mn,x})_{\Omega_{f,x}}| &\leq \sqrt{\varepsilon}\|(F_2-IF_2)_x\|_{\Omega_{f,x}} \sqrt{\varepsilon}\|G_{mn,x}\|_{\Omega_{f,x}} \nonumber\\
    & \leq  C \varepsilon^{1/2} \; h_y^2 \|{F_2}_{xyy}\|_{\Omega_{f,x}} \; 
 \; \|G_{mn}\|_{{_{1,\varepsilon}} } \nonumber  \\
    &\leq C  \frac{\varepsilon^{1/2}}{N^2} \frac{1}{\varepsilon} \|e^{-\frac{\beta}{\sqrt{\varepsilon}}(1+y)}+e^{-\frac{\beta}{\sqrt{\varepsilon}}(1-y)}\|_{L^\infty(\Omega_{f,x})} \sqrt{|\Omega_{f,x}|}\;\sqrt{N \log\left(\frac{1}{\varepsilon}\right)} \nonumber\\
    &\leq C \varepsilon^{1/2} \frac{1}{N^{3/2} \varepsilon} \varepsilon^{1/2} \log (1/\varepsilon)\; \varepsilon^{\frac{3\beta}{2}} \leq C  \frac{\log (1/\varepsilon)}{N^{3/2}}.
\end{align} 
Now, estimate second term of the R.H.S. of  (\ref{F_2}) in the $y$-direction   is $|\varepsilon ((F_2-IF_2)_y,G_{mn,y})_{\Omega_{f,x}}|$. By using Theorem \ref{decompositionu}, Theorem \ref{ajeev}, Lemma \ref{B_1} and furthermore, we obtain
\begin{align}
    |\varepsilon ((F_2-IF_2)_y, G_{mn,y})_{\Omega_{f,x}}| &\leq \sqrt{\varepsilon}\|(F_2-IF_2)_y\|_{\Omega_{f,x}} \sqrt{\varepsilon}\|G_{mn,x}\|_{\Omega_{f,x}} \nonumber\\
    & \leq  C \varepsilon^{1/2} \; h_x^2 \|{F_2}_{yxx}\|_{\Omega_{f,x}} \; 
 \; \|G_{mn}\|_{{_{1,\varepsilon}} } \nonumber  \\
    &\leq C \varepsilon^{1/2} \left(\frac{\varepsilon \log(1/\varepsilon)}{N}\right)^2 \frac{1}{\varepsilon^{1/2}} \|e^{-\frac{\beta}{\sqrt{\varepsilon}}(1+y)}+e^{-\frac{\beta}{\sqrt{\varepsilon}}(1-y)}\|_{L^\infty(\Omega_{f,x})} \nonumber\\
    & \quad \quad\quad\quad\quad\quad\quad\quad\quad\quad\quad\sqrt{|\Omega_{f,x}|}\;\sqrt{N \log\left(\frac{1}{\varepsilon}\right)} \nonumber\\ 
    &\leq C\frac{1}{N^{3/2}} \varepsilon^{2} (\log (1/\varepsilon))^{5/2} 
    \leq C \frac{(\log (1/\varepsilon))^{5/2}}{N^{3/2}}.
\end{align}
 Now, estimate the third term of the R.H.S. of  (\ref{F_2}) in the $x$-direction is
  $|\varepsilon ((F_2-IF_2)_x,G_{mn,x})_{\Omega_{f,y}}|$. By using Theorem \ref{decompositionu}, Theorem \ref{ajeev}, Corollary \ref{B_3} and furthermore, we obtain
\begin{align}
    |\varepsilon ((F_2-IF_2)_x, G_{mn,x})_{\Omega_{f,y}}| &\leq \sqrt{\varepsilon}\|(F_2-IF_2)_x\|_{\Omega_{f,y}} \sqrt{\varepsilon}\|G_{mn,x}\|_{\Omega_{f,y}}  \nonumber\\
    & \leq  C \varepsilon^{1/2} \; h_y^2 \|{F_2}_{xyy}\|_{\Omega_{f,y}} \; 
 \; \|G_{mn}\|_{{_{1,\varepsilon}} }\nonumber  \\
    &\leq C \varepsilon^{1/2} \left(\frac{\sqrt{\varepsilon}\log(1/\varepsilon^{3/2})}{N} \right)^2 \frac{1}{\varepsilon} \|e^{-\frac{\beta}{\sqrt{\varepsilon}}(1+y)}+e^{-\frac{\beta}{\sqrt{\varepsilon}}(1-y)}\|_{L^\infty(\Omega_{f,y})} \nonumber\\
    &\leq C \varepsilon^{1/2} \frac{(\log (1/\varepsilon^{3/2})^{2}}{N} \varepsilon^{1/4} (\log (1/\varepsilon^{3/2})^{1/2} \sqrt{N \log\left(\frac{1}{\varepsilon}\right)}\nonumber\\
    &\leq C \varepsilon^{3/4} \frac{(\log (1/\varepsilon^{3/2}))^{3}}{N^{3/2}} 
     \leq C \frac{(\log (1/\varepsilon^{3/2}))^{3}}{N^{3/2}}.
\end{align}
Now, estimate the third term of the R.H.S. of  (\ref{F_2}) in the $y$-direction is 
  $|\varepsilon ((F_2-IF_2)_y,G_{mn,y})_{\Omega_{f,y}}|$. By using Theorem \ref{decompositionu}, Theorem \ref{ajeev}, Corollary \ref{B_3} and furthermore, we obtain
\begin{align}
    |\varepsilon ((F_2-IF_2)_y, G_{mn,y})_{\Omega_{f,y}}| &\leq \sqrt{\varepsilon}\|(F_2-IF_2)_y\|_{\Omega_{f,y}} \sqrt{\varepsilon}\|G_{mn,y}\|_{\Omega_{f,y}} \nonumber\\
    & \leq  C \varepsilon^{1/2} \; h_x^2 \|{F_2}_{yxx}\|_{\Omega_{f,y}} \; 
 \; \|G_{mn}\|_{{_{1,\varepsilon}} }\nonumber  \\
    &\leq   \frac{C\varepsilon^{1/2}}{N^2} \frac{1}{\varepsilon^{1/2}} \|e^{-\frac{\beta}{\sqrt{\varepsilon}}(1+y)}+e^{-\frac{\beta}{\sqrt{\varepsilon}}(1-y)}\|_{L^\infty(\Omega_{f,y})} \sqrt{|\Omega_{f,y}|}\;\sqrt{N \log\left(\frac{1}{\varepsilon}\right)} \nonumber\\
   & \leq C\frac{1}{N^{3/2}} \varepsilon^{1/4} \log (1/\varepsilon^{3/2}) \leq C \frac{\log (1/\varepsilon^{3/2})}{N^{3/2}}.
\end{align}
 Now, estimate the fourth  term of the R.H.S. of  (\ref{F_2})  in the $x$-direction is
  $|\varepsilon ((F_2-IF_2)_x,G_{mn,x})_{\Omega_{c}}|$. By using Theorem \ref{decompositionu}, Theorem \ref{ajeev}, Lemma \ref{A_1} and furthermore, we obtain
\begin{align}
    |\varepsilon ((F_2-IF_2)_x, G_{mn,x})_{\Omega_{c}}| &\leq \sqrt{\varepsilon}\|(F_2-IF_2)_x\|_{\Omega_{c}} \sqrt{\varepsilon}\|G_{mn,x}\|_{\Omega_{c}}  \nonumber\\
    & \leq  C \varepsilon^{1/2} \; h_y^2 \|{F_2}_{xyy}\|_{\Omega_{c}} \; 
 \; \|G_{mn}\|_{{_{1,\varepsilon}} } \nonumber  \\
    &\leq C  \frac{\sqrt{\varepsilon}}{\varepsilon N^2}  \|e^{-\frac{\beta}{\sqrt{\varepsilon}}(1+y)}+e^{-\frac{\beta}{\sqrt{\varepsilon}}(1-y)}\|_{L^\infty(\Omega_{c})} \sqrt{|\Omega_{c}|}\;\sqrt{N \log\left(\frac{1}{\varepsilon}\right)} \nonumber\\
     &   \leq  \frac{C \varepsilon^{\frac{3\beta}{2}-\frac{1}{2}}}{N^{3/2}}(\log(1/\varepsilon))^{1/2}\leq \frac{C(\log(1/\varepsilon))^{1/2}}{N^{3/2}} \quad \text{when} \; \beta\geq 1/3.
\end{align}
 Now, estimate the fourth  term of the R.H.S. of  (\ref{F_2})  in the $y$-direction is 
  $|\varepsilon ((F_2-IF_2)_y,G_{mn,x})_{\Omega_{c}}|$. By using Theorem \ref{decompositionu}, Theorem \ref{ajeev}, Lemma \ref{A_1} and furthermore, we obtain
\begin{align}
    |\varepsilon ((F_2-IF_2)_y, G_{mn,y})_{\Omega_{c}}| &\leq \sqrt{\varepsilon}\|(F_2-IF_2)_y\|_{\Omega_{c}} \sqrt{\varepsilon}\|G_{mn,y}\|_{\Omega_{c}} \nonumber\\
    & \leq  C \varepsilon^{1/2} \; h_x^2 \|{F_2}_{yxx}\|_{\Omega_{c}} \; 
 \; \|G_{mn}\|_{{_{1,\varepsilon}} }\nonumber  \\
    &\leq C \varepsilon^{1/2} \frac{1}{N^2} \frac{1}{\varepsilon^{1/2}} \sqrt{N \log\left(\frac{1}{\varepsilon}\right)} 
   \leq \frac{C(\log(1/\varepsilon))^{1/2}}{N^{3/2}}.
\end{align}
Applying similar arguments we can show that when source point $(x_m,y_n) \in \Omega_{c}$, then we have
$$|\varepsilon (\nabla(F_2-IF_2),\nabla G_{mn})| \leq C \frac{(\log(1/\varepsilon^{3/2})^{3}}{N}.$$
Hence, when source point $(x_m,y_n) \in \Omega_{c}\cup \Omega_{f,x}$, we have 
$$|\varepsilon (\nabla(F_2-IF_2),\nabla G_{mn})| \leq C \frac{(\log(1/\varepsilon^{3/2})^{3}}{N}.$$
\begin{lemma} \label{diffusion F_12}
     If the source point $(x_m,y_n)\in\Omega_c\cup \Omega_{f,x}$, then for $\alpha\geq2, \; \beta\geq 2/3$ following estimate hold
         $|\varepsilon (\nabla(F_{12}-IF_{12}), \nabla G_{mn})| \leq C  \frac{(\log(1/\varepsilon^{3/2}))^{7/2}}{N}.$
     \end{lemma}
\textbf{Proof.} First we take the source point $(x_m,y_n)\in \Omega_{f,x}$ and decompose the domain in four parts, we have
\begin{align}
    |\varepsilon (\nabla(F_{12}-IF_{12}),\nabla G_{mn})|
    & \leq |\varepsilon (\nabla(F_{12}-IF_{12}),\nabla G_{mn})_{\Omega_{f,xy}}|  
    + |\varepsilon (\nabla(F_{12}-IF_{12}),\nabla G_{mn})_{\Omega_{f,x}}|  \nonumber\\
    &+ |\varepsilon (\nabla(F_{12}-IF_{12}),\nabla G_{mn})_{\Omega_{f,y}}|  
    + |\varepsilon (\nabla(F_{12}-IF_{12}),\nabla G_{mn})_{\Omega_{c}}|   \label{F_12}
\end{align}
Now, estimate  first term of the R.H.S. of (\ref{F_12})  in the $x$-direction   is
  $|\varepsilon ((F_{12}-IF_{12})_x,G_{mn,x})_{\Omega_{f,xy}}|$. By using Theorem \ref{decompositionu}, Theorem \ref{ajeev}, Corollary \ref{A_3}, Corollary \ref{B_3} and furthermore, we obtain
\begin{align}
    |\varepsilon ((F_{12}-IF_{12})_x, G_{mn,x})_{\Omega_{f,xy}}| &\leq \sqrt{\varepsilon}\|(F_{12}-IF_{12})_x\|_{\Omega_{f,xy}} \sqrt{\varepsilon}\|G_{mn,x}\|_{\Omega_{f,xy}} \nonumber\\
    & \leq  C \varepsilon^{1/2} \; h_y^2 \|{F_{12}}_{xyy}\|_{\Omega_{f,xy}} \; 
 \; \|G_{mn}\|_{{_{1,\varepsilon}} }\nonumber  \\
    &\leq C  \left(\frac{\sqrt{\varepsilon}\log(1/\varepsilon^{3/2})}{N} \right)^2  \|e^{-\frac{\alpha x}{\varepsilon}}(e^{-\frac{\beta}{\sqrt{\varepsilon}}(1+y)}+e^{-\frac{\beta}{\sqrt{\varepsilon}}(1-y)})\|_{L^\infty(\Omega_{f,xy})} \nonumber\\
    & \quad \quad  \quad \quad \quad \quad \quad \quad \quad \quad  \frac{\varepsilon^{1/2}}{\varepsilon^2}\sqrt{|\Omega_{f,xy}|}\;\sqrt{N \log\left(\frac{1}{\varepsilon}\right)}\nonumber\\
    &\leq C  \varepsilon^{5/4}\frac{(\log \left(1/\varepsilon^{3/2}\right))^{4}}{N^{3/2}}\leq C \frac{(\log (1/\varepsilon^{3/2}))^{4}}{N^{3/2}}.
\end{align}
Now, estimate first term of  the R.H.S. of (\ref{F_12}) in the $y$-direction   is  
  $|\varepsilon ((F_{12}-IF_{12})_y,G_{mn,y})_{\Omega_{f,xy}}|$. By using Theorem \ref{decompositionu}, Theorem \ref{ajeev}, Corollary \ref{A_3}, Corollary \ref{B_3} and furthermore, we obtain
  \begin{align}
    |\varepsilon ((F_{12}-IF_{12})_y, G_{mn,y})_{\Omega_{f,xy}}| &\leq \sqrt{\varepsilon}\|(F_{12}-IF_{12})_y\|_{\Omega_{f,xy}} \sqrt{\varepsilon}\|G_{mn,y}\|_{\Omega_{f,xy}} \nonumber\\
    & \leq  C \varepsilon^{1/2} \; h_x^2 \|{F_{12}}_{yxx}\|_{\Omega_{f,xy}} \; 
 \; \|G_{mn}\|_{{_{1,\varepsilon}} }\nonumber  \\
    &\leq C \varepsilon^{1/2}  \left(\frac{\varepsilon\log(1/\varepsilon)}{N} \right)^2\frac{\varepsilon^{3/4}}{\varepsilon^{2}}  \log(1/\varepsilon^{3/2})\sqrt{N \log\left(\frac{1}{\varepsilon}\right)} \nonumber\\
    &\leq C \varepsilon^{5/4}\frac{ (\log (1/\varepsilon^{3/2}))^{7/2}}{N^{3/2}} \leq C\frac{ (\log (1/\varepsilon^{3/2}))^{7/2}}{N^{3/2}} .
\end{align}
 Now, estimate second term of the R.H.S. of (\ref{F_12}) in the $x$-direction   is
  $|\varepsilon ((F_{12}-IF_{12})_x,G_{mn,x})_{\Omega_{f,x}}|$. By using Theorem \ref{decompositionu}, Theorem \ref{ajeev}, Corollary \ref{A_3}, Lemma \ref{B_1} and furthermore, we obtain
\begin{align}
    |\varepsilon ((F_{12}-IF_{12})_x, G_{mn,x})_{\Omega_{f,x}}| &\leq \sqrt{\varepsilon}\|(F_{12}-IF_{12})_x\|_{\Omega_{f,x}} \sqrt{\varepsilon}\|G_{mn,x}\|_{\Omega_{f,x}} \nonumber\\
    & \leq  C \varepsilon^{1/2} \; h_y^2 \|{F_{12}}_{xyy}\|_{\Omega_{f,x}} \; 
 \; \|G_{mn}\|_{{_{1,\varepsilon}} }\nonumber  \\
    &\leq  \frac{C\varepsilon^{1/2}}{\varepsilon^2N^2} \|e^{-\frac{\alpha x}{\varepsilon}}(e^{-\frac{\beta}{\sqrt{\varepsilon}}(1+y)}+e^{-\frac{\beta}{\sqrt{\varepsilon}}(1-y)})\|_{L^\infty(\Omega_{f,x})}  \;\|G_{mn}\|_{{_{1,\varepsilon}} } \nonumber\\
   &\leq  \frac{C(\varepsilon\log (\frac{1}{\varepsilon}))^{1/2}}{\varepsilon^{3/2} N} \varepsilon^{\frac{3}{2}\beta}
     \sqrt{N \log\left(\frac{1}{\varepsilon}\right)}\nonumber\\
     & \leq   C  \frac{\log (1/\varepsilon)}{ N} \varepsilon^{\frac{3}{2}\beta-1}\leq  C  \frac{(\log (1/\varepsilon))^{1/2}}{ N} \quad \text{when} \; \beta\geq 2/3.
\end{align}
Now, estimate second term of  the R.H.S. of (\ref{F_12}) in the $y$-direction   is  
  $|\varepsilon ((F_{12}-IF_{12})_y,G_{mn,y})_{\Omega_{f,x}}|$. By using Theorem \ref{decompositionu}, Theorem \ref{ajeev}, Corollary \ref{A_3}, Lemma \ref{B_1} and furthermore, we obtain
 \begin{align}
    |\varepsilon ((F_{12}-IF_{12})_y, G_{mn,y})_{\Omega_{f,x}}| &\leq \sqrt{\varepsilon}\|(F_{12}-IF_{12})_y\|_{\Omega_{f,x}} \sqrt{\varepsilon}\|G_{mn,y}\|_{\Omega_{f,x}} \nonumber\\
    & \leq  C \varepsilon^{1/2} \; h_x^2 \|{F_{12}}_{yxx}\|_{\Omega_{f,x}} \; 
 \; \|G_{mn}\|_{{_{1,\varepsilon}} }\nonumber  \\
    &\leq C   \left(\frac{\varepsilon\log(1/\varepsilon)}{N} \right)^2\frac{\varepsilon^{1/2}}{\varepsilon^{5/2}} \|e^{-\frac{\alpha x}{\varepsilon}}(e^{-\frac{\beta}{\sqrt{\varepsilon}}(1+y)}+e^{-\frac{\beta}{\sqrt{\varepsilon}}(1-y)})\|_{L^\infty(\Omega_{f,x})} \nonumber\\
    & \quad  \quad\quad\quad\quad\quad\quad\quad\quad\quad\sqrt{|\Omega_{f,x}|}\;\sqrt{N \log\left(\frac{1}{\varepsilon}\right)}\nonumber\\
   & \leq C\frac{1}{N^{3/2}} \varepsilon^{1/2} (\log (1/\varepsilon))^{3} \varepsilon^{\frac{3\beta}{2}} \leq C \frac{(\log (1/\varepsilon))^{3}}{N^{3/2}}.
\end{align}
Now, estimate  third term of the R.H.S. of (\ref{F_12}) in the $x$-direction   is 
  $|\varepsilon ((F_{12}-IF_{12})_x,G_{mn,x})_{\Omega_{f,y}}|$. By using Theorem \ref{decompositionu}, Theorem \ref{ajeev}, Lemma \ref{A_1}, Corollary \ref{B_3} and furthermore, we obtain
\begin{align}
    |\varepsilon ((F_{12}-IF_{12})_x, G_{mn,x})_{\Omega_{f,y}}| &\leq \sqrt{\varepsilon}\|(F_{12}-IF_{12})_x\|_{\Omega_{f,y}} \sqrt{\varepsilon}\|G_{mn,x}\|_{\Omega_{f,y}} \nonumber\\
    & \leq  C \varepsilon^{1/2} \; h_y^2 \|{F_{12}}_{xyy}\|_{\Omega_{f,y}} \; 
 \; \|G_{mn}\|_{{_{1,\varepsilon}} }\nonumber  \\
    &\leq C \left(\frac{\sqrt{\varepsilon}\log(1/\varepsilon^{3/2})}{N} \right)^2 \|e^{-\frac{\alpha x}{\varepsilon}}(e^{-\frac{\beta}{\sqrt{\varepsilon}}(1+y)}+e^{-\frac{\beta}{\sqrt{\varepsilon}}(1-y)})\|_{L^\infty(\Omega_{f,y})} \nonumber\\
    & \quad  \quad \quad  \quad \quad \quad \quad \quad \quad\frac{\varepsilon^{1/2}}{\varepsilon^2}\sqrt{|\Omega_{f,y}|}\;\sqrt{N \log\left(\frac{1}{\varepsilon}\right)} \nonumber\\
    &\leq C  \frac{(\log (1/\varepsilon^{3/2}))^{3}}{\varepsilon^{1/4} N^{3/2}} \varepsilon^\alpha  \leq C \frac{(\log (1/\varepsilon^{3/2}))^{3}}{ N^{3/2}} \quad \text{when}\; \alpha\geq1/4.
\end{align}
Now, estimate  third term of   the R.H.S. of (\ref{F_12}) in the $y$-direction   is 
  $|\varepsilon ((F_{12}-IF_{12})_y,G_{mn,y})_{\Omega_{f,y}}|$. By using Theorem \ref{decompositionu}, Theorem \ref{ajeev}, Lemma \ref{A_1}, Corollary \ref{B_3} and furthermore, we obtain
\begin{align}
    |\varepsilon ((F_{12}-IF_{12})_y, G_{mn,y})_{\Omega_{f,y}}| &\leq \sqrt{\varepsilon}\|(F_{12}-IF_{12})_y\|_{\Omega_{f,y}} \sqrt{\varepsilon}\|G_{mn,y}\|_{\Omega_{f,y}} \nonumber\\
    & \leq  C \varepsilon^{1/2} \; h_x^2 \|{F_{12}}_{yxx}\|_{\Omega_{f,y}} \; 
 \; \|G_{mn}\|_{{_{1,\varepsilon}} } \nonumber  \\
    &\leq   \frac{C\sqrt{\varepsilon}}{N^2\varepsilon^{\frac{5}{2}}}\|e^{-\frac{\alpha x}{\varepsilon}}(e^{-\frac{\beta}{\sqrt{\varepsilon}}(1+y)}+e^{-\frac{\beta}{\sqrt{\varepsilon}}(1-y)})\|_{L^\infty(\Omega_{f,y})}  \sqrt{|\Omega_{f,y}|}\sqrt{N\log\frac{1}{\varepsilon}} \nonumber\\
    &\leq C  \frac{\log (1/\varepsilon^{3/2})}{\varepsilon^{7/4} N^{3/2}} \varepsilon^\alpha 
 \leq C \frac{\log (1/\varepsilon^{3/2})}{ N^{3/2}} \quad \text{when}\; \alpha\geq7/4.
\end{align}
Now, estimate  fourth  term of the R.H.S. of  (\ref{F_12}) in the $x$-direction   is  
  $|\varepsilon ((F_{12}-IF_{12})_x,G_{mn,x})_{\Omega_{c}}|$. By using Theorem \ref{decompositionu}, Theorem \ref{ajeev}, Lemma \ref{A_1}, Lemma \ref{B_1} and furthermore, we obtain
\begin{align}
    |\varepsilon ((F_{12}-IF_{12})_x, G_{mn,x})_{\Omega_{c}}| &\leq \sqrt{\varepsilon}\|(F_{12}-IF_{12})_x\|_{\Omega_{c}} \sqrt{\varepsilon}\|G_{mn,x}\|_{\Omega_{c}} \nonumber\\
    & \leq  C \varepsilon^{1/2} \; h_y^2 \|{F_{12}}_{xyy}\|_{\Omega_{c}} \; 
 \; \|G_{mn}\|_{{_{1,\varepsilon}} }\nonumber  \\
    & \frac{\leq C \varepsilon^{1/2}}{N^2}\frac{1}{\varepsilon^2} \|e^{-\frac{\alpha x}{\varepsilon}}(e^{-\frac{\beta}{\sqrt{\varepsilon}}(1+y)}+e^{-\frac{\beta}{\sqrt{\varepsilon}}(1-y)})\|_{L^\infty(\Omega_{c})} \sqrt{|\Omega_{c}|}\sqrt{N\log\frac{1}{\varepsilon}} \nonumber\\
    &\leq   \frac{C}{\varepsilon^{3/2} N^{3/2}}\sqrt{\log\frac{1}{\varepsilon}} \varepsilon^{\alpha+\frac{3\beta}{2}}
    \leq  \frac{C(\log(1/\varepsilon))^{1/2}}{ N^{3/2}} \quad \text{when}\; \alpha+\frac{3\beta}{2}\geq3/2.
\end{align}
Now, estimate fourth  term of  the R.H.S. of (\ref{F_12}) in the $y$-direction   is  
  $|\varepsilon ((F_{12}-IF_{12})_y,G_{mn,y})_{\Omega_{c}}|$. By using Theorem \ref{decompositionu}, Theorem \ref{ajeev}, Lemma \ref{A_1}, Lemma \ref{B_1} and furthermore, we obtain
 \begin{align}
    |\varepsilon ((F_{12}-IF_{12})_y, G_{mn,y})_{\Omega_{c}}| &\leq \sqrt{\varepsilon}\|(F_{12}-IF_{12})_y\|_{\Omega_{c}} \sqrt{\varepsilon}\|G_{mn,y}\|_{\Omega_{c}} \nonumber\\
    & \leq  C \varepsilon^{1/2} \; h_x^2 \|{F_{12}}_{yxx}\|_{\Omega_{c}} \; 
 \; \|G_{mn}\|_{{_{1,\varepsilon}} }\nonumber  \\
    &\leq C  \frac{\varepsilon^{1/2}}{\varepsilon^{\frac{5}{2}}N^2}\|e^{-\frac{\alpha x}{\varepsilon}}(e^{-\frac{\beta}{\sqrt{\varepsilon}}(1+y)}+e^{-\frac{\beta}{\sqrt{\varepsilon}}(1-y)})\|_{L^\infty(\Omega_{c})} \sqrt{|\Omega_{c}|}\sqrt{N\log\frac{1}{\varepsilon}} \nonumber\\
    &\leq C  \frac{(\log(1/\varepsilon))^{1/2}}{ N^{3/2}} \varepsilon^{\alpha+\frac{3\beta}{2}-2}  \leq  \frac{C(\log(1/\varepsilon))^{1/2}}{ N^{3/2}} \quad \text{when}\; \alpha+\frac{3\beta}{2}\geq2.
\end{align}
Applying similar arguments we can show that when source point $(x_m,y_n) \in \Omega_{c}$, then we have
$$|\varepsilon (\nabla(F_{12}-IF_{12}),\nabla G_{mn})| \leq C \frac{(\log(1/\varepsilon^{3/2})^{7/2}}{N}.$$
Hence, when source point $(x_m,y_n) \in \Omega_{c}\cup \Omega_{f,x}$, we have
$$|\varepsilon (\nabla(F_{12}-IF_{12}),\nabla G_{mn})| \leq C \frac{(\log(1/\varepsilon^{3/2})^{7/2}}{N}.$$
\section{Estimates related to convection and reaction terms}\label{section5}
In this section, we aim to bound the convection and reaction portion of the bilinear form B(.,.) with its arguments parts of interpolation error $u-Iu$ and the discrete Green's function $G_{mn}.$ For the convection term, while we get $\varepsilon$-independent estimate for $(x_m,y_n)\in \Omega_{c}\cup\Omega_{f,x}$, the convergence behavior differs. For  $(x_m,y_n)\in \Omega_{c}$, we could not determine a specific order of convergence. For $(x_m,y_n)\in\Omega_{f,x}$, we obtained a convergence order of $N^{-1/2}$ up to a logarithmic factor. And for the reaction term, we got a convergence order of $N^{-1}$ up to a logarithmic factor when $(x_m,y_n)\in \Omega_{c}\cup\Omega_{f,x}$. We begin our section with the following auxiliary results.
\begin{theorem}{\cite{Sun94}} \label{interpolation}
        Assume that $Iu$ represents the bilinear interpolation of $u$ on the Shishkin-type mesh. The interpolation error thus satisfies
\begin{center}
    $ \|u-Iu\|_{L^\infty(\Omega_f)}\leq CN^{-2} \log^2(1/\varepsilon^{3/2}), \quad\quad \|u-Iu\|_{L^\infty(\Omega_c)}\leq CN^{-2}.$
\end{center}
\end{theorem}
\begin{lemma} \label{convection term f,x}
     If the source point $(x_m,y_n)\in\Omega_c\cup \Omega_{f,x}$, then for $  \beta\geq 2/3$ following estimate hold
 $$|(Iu-u,b_1G_{mn,x})_{\Omega_{f,x}}| \leq C \frac{(\log(1/\varepsilon))^{2}}{N}.$$
\end{lemma}
\textbf{Proof.} We begin by decomposing the solution $u$ into its regular and singular parts as $u=R+F_1+F_2+F_{12}$.
Consider the source point $(x_m,y_n)\in \Omega_{f,x}$. Hence,
  \begin{align}
    |(Iu-u,b_1G_{mn,x})_{\Omega_{f,x}}|
     \leq & |(IR-R,b_1G_{mn,x})_{\Omega_{f,x}}|
     +  |(IF_1-F_1,b_1G_{mn,x})_{\Omega_{f,x}}| \nonumber\\
    &+  |(IF_2-F_2,b_1G_{mn,x})_{\Omega_{f,x}}| 
    + |(IF_{12}-F_{12},b_1G_{mn,x})_{\Omega_{f,x}}|.  \label{Convection f,x}
\end{align}
Consider the first term of the R.H.S. of (\ref{Convection f,x}) is $|(IR-R,b_1G_{mn,x})_{\Omega_{f,x}}|$. By using Theorem \ref{interpolation} and furthermore, we obtain
\begin{align}
    |(IR-R,b_1G_{mn,x})_{\Omega_{f,x}}| & \leq \|IR-R\|_{\Omega_{f,x}} \|b_1 G_{mn,x}\|_{\Omega_{f,x}} \nonumber\\
    & \leq C \|IR-R\|_{L^{\infty}( \Omega_{f,x})} \sqrt{|\Omega_{f,x}|} \; \| G_{mn,x}\|_{\Omega_{f,x}} \nonumber\\
     &\leq \frac{C}{N^2} \log^{1/2}(1/\varepsilon)\;\| G_{mn}\|_{1,\varepsilon} \leq \frac{C}{N^2} \log^{1/2}(1/\varepsilon)\sqrt{N \log\left(\frac{1}{\varepsilon}\right)} \nonumber\\
     & \leq \frac{C}{N^{3/2}} \log(1/\varepsilon).
\end{align}
Consider the second term of the R.H.S. of (\ref{Convection f,x}) is $|(IF_1-F_1,b_1G_{mn,x})_{\Omega_{f,x}}|$. By using Theorem \ref{decompositionu}, Theorem \ref{ajeev}, Corollary \ref{A_3} and furthermore, we obtain
\begin{align}
    |(IF_1-F_1,b_1G_{mn,x})_{\Omega_{f,x}}| & \leq \|IF_1-F_1\|_{\Omega_{f,x}} \|b_1 G_{mn,x}\|_{\Omega_{f,x}} \nonumber\\
    & \leq C \|IF_1-F_1\|_{L^{\infty}( \Omega_{f,x})} \sqrt{|\Omega_{f,x}|} \; \| G_{mn,x}\|_{\Omega_{f,x}} \nonumber\\
      &\leq C\|IF_1-F_1\|_{L^{\infty}( \Omega_{f,x})} \log^{1/2}(1/\varepsilon)\;\| G_{mn}\|_{1,\varepsilon} \nonumber\\
      &\leq C\sqrt{N \log\left(\frac{1}{\varepsilon}\right)} [h_x^2\|F_{1,xx}\|_{L^{\infty}( \Omega_{f,x})}+h_x h_y\|F_{1,xy}\|_{L^{\infty}( \Omega_{f,x})}\nonumber\\
      &\quad\quad\quad\quad\quad\quad\quad +h_y^2\|F_{1,yy}\|_{L^{\infty}( \Omega_{f,x})}] \nonumber\\
      &\leq C\sqrt{N \log\left(\frac{1}{\varepsilon}\right)}\left[\frac{\varepsilon^2\log^{2}(1/\varepsilon)}{\varepsilon^2N^2} +\frac{\varepsilon\log^{2}(1/\varepsilon)}{N^2\varepsilon} +\frac{1}{N^2}\right]\nonumber\\
      &\leq C \frac{(\log(1/\varepsilon))^{5/2}}{N^{3/2}}.
\end{align}
Next consider the third term of the R.H.S. of (\ref{Convection f,x}) is $|(IF_2-F_2,b_1G_{mn,x})_{\Omega_{f,x}}|$. By using Theorem \ref{decompositionu}, Theorem \ref{ajeev}, Lemma \ref{B_1} and furthermore, we obtain 
\begin{align}
    |(IF_2-F_2,b_1G_{mn,x})_{\Omega_{f,x}}| & \leq \|IF_2-F_2\|_{\Omega_{f,x}} \|b_1 G_{mn,x}\|_{\Omega_{f,x}} \nonumber\\
    & \leq C \|IF_2-F_2\|_{L^{\infty}( \Omega_{f,x})} \sqrt{|\Omega_{f,x}|} \; \| G_{mn,x}\|_{\Omega_{f,x}} \nonumber\\
      &\leq C\|IF_2-F_2\|_{L^{\infty}( \Omega_{f,x})} \log^{1/2}(1/\varepsilon)\;\| G_{mn}\|_{1,\varepsilon} \nonumber\\
      &\leq CN [h_x^2\|F_{2,xx}\|_{L^{\infty}( \Omega_{f,x})}+h_x h_y\|F_{2,xy}\|_{L^{\infty}( \Omega_{f,x})}+h_y^2\|F_{2,yy}\|_{L^{\infty}( \Omega_{f,x})}]  \nonumber\\
      &\leq C\sqrt{N \log\left(\frac{1}{\varepsilon}\right)} \left[\frac{\varepsilon^2\log^{2}(1/\varepsilon)}{N^2} +\frac{\varepsilon\log^{2}(1/\varepsilon)}{N^2} \frac{1}{\varepsilon^{1/2}} \right]\nonumber\\
      & \quad\quad\quad\quad\quad + C\sqrt{N \log\left(\frac{1}{\varepsilon}\right)} \left[\frac{1}{N^2\varepsilon}\|e^{-\frac{\beta}{\sqrt{\varepsilon}}(1+y)}+e^{-\frac{\beta}{\sqrt{\varepsilon}}(1-y)}\|_{L^\infty(\Omega_{f,x})} \right]\nonumber\\
      &\leq  C\sqrt{N \log\left(\frac{1}{\varepsilon}\right)} \left[\frac{\varepsilon^2\log^{2}(1/\varepsilon)}{N^2} +\frac{\varepsilon\log^{2}(1/\varepsilon)}{N^2} \frac{1}{\varepsilon^{1/2}} + \frac{\varepsilon^{\frac{3\beta}{2}-1}}{N^2}\right]\nonumber\\
      &\leq  C \frac{(\log(1/\varepsilon))^{5/2}}{N^{3/2}}.
\end{align}
Next consider the fourth term of the R.H.S. of (\ref{Convection f,x}) is $|(IF_{12}-F_{12},b_1G_{mn,x})_{\Omega_{f,x}}|$. By using Theorem \ref{decompositionu}, Theorem \ref{ajeev}, Corollary \ref{A_3}, Lemma \ref{B_1} and furthermore, we obtain
\begin{align}
    |(IF_{12}-F_{12},b_1G_{mn,x})_{\Omega_{f,x}}| & \leq \|IF_{12}-F_{12}\|_{\Omega_{f,x}} \|b_1 G_{mn,x}\|_{\Omega_{f,x}} \nonumber\\
    & \leq C \|IF_{12}-F_{12}\|_{L^{\infty}( \Omega_{f,x})} \sqrt{|\Omega_{f,x}|} \; \| G_{mn,x}\|_{\Omega_{f,x}} \nonumber\\
      &\leq C\|IF_{12}-F_{12}\|_{L^{\infty}( \Omega_{f,x})} \log^{1/2}(1/\varepsilon)\;\| G_{mn}\|_{1,\varepsilon} \nonumber\\
      &\leq C [h_x^2\|F_{12,xx}\|_{L^{\infty}( \Omega_{f,x})}+h_x h_y\|F_{12,xy}\|_{L^{\infty}( \Omega_{f,x})}  \nonumber\\
       &  \quad\quad\quad\quad\quad+h_y^2\|F_{12,yy}\|_{L^{\infty}( \Omega_{f,x})}] \sqrt{N \log\left(\frac{1}{\varepsilon}\right)}\nonumber\\
      &\leq C \sqrt{N \log\left(\frac{1}{\varepsilon}\right)}\left[h_x^2 +\frac{h_xh_y}{\varepsilon^{3/2}} +\frac{h_y^2}{\varepsilon} \right]\varepsilon^{\frac{3\beta}{2}}\nonumber\\
      &\leq  C\sqrt{N \log\left(\frac{1}{\varepsilon}\right)}\left[\frac{\varepsilon^2\log^{2}(1/\varepsilon)\varepsilon^{\frac{3\beta}{2}}}{N^2} +\frac{\varepsilon\log(1/\varepsilon)}{N^2} \frac{\varepsilon^{\frac{3\beta}{2}}}{\varepsilon^{3/2}} + \frac{\varepsilon^{\frac{3\beta}{2}-1}}{N^2}\right]\nonumber\\
      &\leq  C \frac{(\log(1/\varepsilon))^{5/2}}{N^{3/2}}. 
\end{align}
Hence when source point $(x_m,y_n) \in \Omega_{f,x}$, then we have 
$$|(Iu-u,b_1G_{mn,x})_{\Omega_{f,x}}| \leq C \frac{(\log(1/\varepsilon))^{5/2}}{N^{3/2}}.$$
Similar arguments we can show that when source point $(x_m,y_n) \in \Omega_{c}$, then we have
$$|(Iu-u,b_1G_{mn,x})_{\Omega_{f,x}}| \leq C \frac{(\log(1/\varepsilon))^2}{N}.$$
Hence, when source point $(x_m,y_n) \in \Omega_{c}\cup \Omega_{f,x}$, we have 
$$|(Iu-u,b_1G_{mn,x})_{\Omega_{f,x}}| \leq C \frac{(\log(1/\varepsilon))^{2}}{N}.$$
\begin{lemma} \label{convection term f,y}
   If the source point $(x_m,y_n)\in\Omega_c\cup \Omega_{f,x}$, then for $  \alpha\geq2$ following estimate hold
$$|(Iu-u,b_1G_{mn,x})_{\Omega_{f,y}}| \leq C \frac{(\log(1/{\varepsilon}))^{1/2}}{N}.$$
\end{lemma}
\textbf{Proof.} First we take the source point $(x_m,y_n)\in \Omega_{f,x}$ and consider $ |(Iu-u,b_1G_{mn,x})_{\Omega_{f,y}}|$ and decompose $u$ as $u=R+F_1+F_2+F_{12}.$ Hence,
  \begin{align}
    |(Iu-u,b_1G_{mn,x})_{\Omega_{f,y}}|
     \leq  & |(IR-R,b_1G_{mn,x})_{\Omega_{f,y}}|  +  |(IF_1-F_1,b_1G_{mn,x})_{\Omega_{f,y}}| \nonumber\\
    +  &|(IF_2-F_2,b_1G_{mn,x})_{\Omega_{f,y}}| +|(IF_{12}-F_{12},b_1G_{mn,x})_{\Omega_{f,y}}|  \label{Convection f,y}
\end{align}
Consider the second term of the R.H.S. of (\ref{Convection f,y}) is $|(IF_1-F_1,b_1G_{mn,x})_{\Omega_{f,y}}|$. By using inverse inequality, Theorem \ref{decompositionu}, Theorem \ref{ajeev}, Lemma \ref{A_1} and furthermore, we obtain
\begin{align}
    |(IF_1-F_1,b_1G_{mn,x})_{\Omega_{f,y}}| & \leq \|IF_1-F_1\|_{\Omega_{f,y}} \|b_1 G_{mn,x}\|_{\Omega_{f,y}} \nonumber\\
    & \leq C \|IF_1-F_1\|_{L^{\infty}( \Omega_{f,y})} \sqrt{|\Omega_{f,y}|} \; \| G_{mn,x}\|_{\Omega_{f,y}} \nonumber\\
      &\leq C\|IF_1-F_1\|_{L^{\infty}( \Omega_{f,y})} \log^{1/2}(1/\varepsilon^{3/2})\varepsilon^{1/4}\;N\| G_{mn}\|_{\Omega_{f,y}} \nonumber\\
      &\leq CN\varepsilon^{1/4}  [h_x^2\|F_{1,xx}\|_{L^{\infty}( \Omega_{f,y})}+h_x h_y\|F_{1,xy}\|_{L^{\infty}( \Omega_{f,y})} \nonumber\\
      &\quad+h_y^2\|F_{1,yy}\|_{L^{\infty}( \Omega_{f,y})}] (N \log(1/{\varepsilon}))^{1/2} \nonumber\\
      &\leq CN(N \log(1/{\varepsilon}))^{1/2}\varepsilon^{1/4} \left[\frac{1}{N^2\varepsilon^2} +\frac{h_x h_y}{\varepsilon} +h_y^2\right]\|e^{\frac{-\alpha x}{\varepsilon}}\|_{L^{\infty}( \Omega_{f,y})} \nonumber\\
      &\leq \frac{C\varepsilon^{1/4}( \log(1/{\varepsilon}))^{1/2}}{\varepsilon^2 N^{1/2}}\varepsilon^{\alpha}\leq \frac{C(\log(1/{\varepsilon}))^{1/2}}{N}\quad \text{when} ~ \alpha \geq2.
\end{align}
Next consider the third term of the R.H.S. of (\ref{Convection f,y}) is $|(IF_2-F_2,b_1G_{mn,x})_{\Omega_{f,y}}|$. By using inverse inequality, Theorem \ref{decompositionu}, Theorem \ref{ajeev}, Corollary \ref{B_3} and furthermore, we obtain
\begin{align}
    |(IF_2-F_2,b_1G_{mn,x})_{\Omega_{f,y}}| & \leq \|IF_2-F_2\|_{\Omega_{f,y}} \|b_1 G_{mn,x}\|_{\Omega_{f,y}} \nonumber\\
    & \leq C \|IF_2-F_2\|_{L^{\infty}( \Omega_{f,y})} \sqrt{|\Omega_{f,y}|} \; \| G_{mn,x}\|_{\Omega_{f,y}} \nonumber\\
      &\leq C\|IF_2-F_2\|_{L^{\infty}( \Omega_{f,y})} \log^{1/2}(1/\varepsilon^{3/2})\;\varepsilon^{1/4}N\| G_{mn}\|_{\Omega_{f,y}} \nonumber\\
      &\leq CN(N \log(1/{\varepsilon}))^{1/2}\varepsilon^{1/4} [h_x^2\|F_{2,xx}\|_{L^{\infty}( \Omega_{f,y})}+h_x h_y\|F_{2,xy}\|_{L^{\infty}( \Omega_{f,y})}\nonumber\\
      &\quad+h_y^2\|F_{2,yy}\|_{L^{\infty}( \Omega_{f,y})}]  \nonumber\\
      &\leq CN(N \log(1/{\varepsilon}))^{1/2}\varepsilon^{1/4}\left[\frac{1}{N^2} +\frac{\sqrt{}{\varepsilon}\log(1/\varepsilon^{3/2})}{\varepsilon^{1/2}N^2}  + \frac{\varepsilon\log^2{(1/\varepsilon^{3/2})}}{\varepsilon N^2}  \right]\nonumber\\
      & \leq \frac{C\varepsilon^{1/4} (\log(1/{\varepsilon}))^{1/2}}{N^{1/2}} \leq \frac{C (\log(1/{\varepsilon}))^{1/2}}{N}.
\end{align}
Next consider the fourth term of the R.H.S. of (\ref{Convection f,y}) is $|(IF_{12}-F_{12},b_1G_{mn,x})_{\Omega_{f,y}}|$. By using Theorem \ref{decompositionu}, Theorem \ref{ajeev}, Lemma \ref{A_1}, Corollary \ref{B_3} and furthermore, we obtain
\begin{align}
    |(IF_{12}-F_{12},b_1G_{mn,x})_{\Omega_{f,y}}| & \leq \|IF_{12}-F_{12}\|_{\Omega_{f,y}} \|b_1 G_{mn,x}\|_{\Omega_{f,y}} \nonumber\\
    & \leq C \|IF_{12}-F_{12}\|_{L^{\infty}( \Omega_{f,y})} \sqrt{|\Omega_{f,y}|} \; \| G_{mn,x}\|_{\Omega_{f,y}} \nonumber\\
      &\leq C\|IF_{12}-F_{12}\|_{L^{\infty}( \Omega_{f,y})} \varepsilon^{1/4}\log^{1/2}(1/\varepsilon^{3/2})\;N\| G_{mn}\|_{\Omega_{f,y}} \nonumber\\
      &\leq CN \varepsilon^{1/4} [h_x^2\|F_{12,xx}\|_{L^{\infty}( \Omega_{f,y})}+h_x h_y\|F_{12,xy}\|_{L^{\infty}( \Omega_{f,y})}  \nonumber\\
       &  \quad\quad\quad\quad\quad+h_y^2\|F_{12,yy}\|_{L^{\infty}( \Omega_{f,y})}] (N \log(1/{\varepsilon}))^{1/2}\nonumber\\
      &\leq CN(N \log(1/{\varepsilon}))^{1/2} \left[{\frac{\varepsilon^{1/4}}{\varepsilon^{2}N^2}} +\frac{h_y h_x\varepsilon^{1/4}}{\varepsilon^{3/2}} +\frac{\varepsilon\log^2(1/\varepsilon^{3/2})}{N^2\varepsilon} \right] \varepsilon^{\alpha}\nonumber\\
      &\leq  CN(N \log(1/{\varepsilon}))^{1/2} \varepsilon^{1/4}\left[\frac{\varepsilon^{\alpha}}{\varepsilon^{2}N^2}+\frac{\sqrt{\varepsilon}\log(1/\varepsilon^{3/2})}{N^2} \frac{\varepsilon^\alpha}{\varepsilon^{3/2}} + \frac{\varepsilon^{\alpha}}{N^2}\right]\nonumber\\
      &\leq \frac{C\varepsilon^{1/4} (\log(1/{\varepsilon}))^{1/2}}{N^{1/2}}\leq  C \frac{( \log(1/{\varepsilon}))^{1/2}}{N}. 
\end{align}
Consider the first term of the R.H.S. of (\ref{Convection f,y}) is $|(IR-R,b_1G_{mn,x})_{\Omega_{f,y}}|$. By using Theorem \ref{decompositionu}, Theorem \ref{ajeev}, Corollary \ref{A_3}, Corollary \ref{B_3} and furthermore, we obtain
\begin{align}
    |(IR-R,b_1G_{mn,x})_{\Omega_{f,y}}| & \leq \|IR-R\|_{\Omega_{f,y}} \|b_1 G_{mn,x}\|_{\Omega_{f,y}} \nonumber\\
    & \leq C \|IR-R\|_{L^{\infty}( \Omega_{f,y})} \sqrt{|\Omega_{f,y}|} \; \| G_{mn,x}\|_{\Omega_{f,y}} \nonumber\\
     &\leq C\|IR-R\|_{L^{\infty}( \Omega_{f,y})} \log(1/\varepsilon^{3/2})\varepsilon^{1/4}\;\| G_{mn}\|_{1,\varepsilon} \nonumber\\
      &\leq \frac{C\|G_{mn}\|_{1,\varepsilon}}{\varepsilon^{1/4}} [h_x^2\|R_{xx}\|_{L^{\infty}( \Omega_{f,y})}+h_x h_y\|R_{xy}\|_{L^{\infty}( \Omega_{f,y})}+h_y^2\|R_{yy}\|_{L^{\infty}( \Omega_{f,y})}]  \nonumber\\
      &\leq \frac{C(N \log(1/{\varepsilon}))^{1/2}}{\varepsilon^{1/4}} \left[\frac{\varepsilon^2 \log^2(1/\varepsilon)}{N^2} +\frac{\varepsilon\sqrt{\varepsilon}\log^2(1/\varepsilon^{3/2})}{ N^2} +\frac{\varepsilon \log^2(1/\varepsilon^{3/2})}{N^2} \right] \nonumber\\
       &\leq C \frac{(\log({1/\varepsilon^{3/2}))^{5/2}}}{N^{3/2}}.
\end{align}
Hence, when source point $(x_m,y_n) \in \Omega_{f,x}$, then 
$$|(Iu-u,b_1G_{mn,x})_{\Omega_{f,y}}| \leq C \frac{( \log(1/{\varepsilon}))^{1/2}}{N}.$$
Similar arguments we can show that when source point $(x_m,y_n) \in \Omega_{c}$, then we have 
$$|(Iu-u,b_1G_{mn,x})_{\Omega_{f,y}}| \leq C \frac{\log^{2}(1/\varepsilon)}{N}.$$
Hence, when source point $(x_m,y_n) \in \Omega_{c}\cup \Omega_{f,x}$, we have 
$$|(Iu-u,b_1G_{mn,x})_{\Omega_{f,y}}| \leq C \frac{( \log(1/{\varepsilon}))^{1/2}}{N}.$$
\begin{lemma}\label{convection term f,xy}
    If the source point $(x_m,y_n)$ lies in $\Omega_c\cup \Omega_{f,x}$, then following estimate hold
$$|(Iu-u,b_1G_{mn,x})_{\Omega_{f,xy}}| \leq C \frac{(\log({1/\varepsilon^{3/2}))^{5/2}}}{N}.$$
\end{lemma}
\textbf{Proof.} First we take the source point $(x_m,y_n)\in \Omega_{f,x}$ and consider $|(Iu-u,b_1G_{mn,x})_{\Omega_{f,xy}}|$ and decompose $u$ as $u=R+F_1+F_2+F_{12}.$ Hence,
  \begin{align}
    |(Iu-u,b_1G_{mn,x})_{\Omega_{f,xy}}|
     \leq & |(IR-R,b_1G_{mn,x})_{\Omega_{f,xy}}| 
     +  |(IF_1-F_1,b_1G_{mn,x})_{\Omega_{f,xy}}| \nonumber\\
    &+  |(IF_2-F_2,b_1G_{mn,x})_{\Omega_{f,xy}}| 
    + |(IF_{12}-F_{12},b_1G_{mn,x})_{\Omega_{f,xy}}|.  \label{Convection f,xy}
\end{align}
Consider the first term of the R.H.S. of (\ref{Convection f,xy}) is $|(IR-R,b_1G_{mn,x})_{\Omega_{f,xy}}|$. By using Theorem \ref{decompositionu}, Theorem \ref{ajeev}, Corollary \ref{A_3}, Corollary \ref{B_3} and furthermore, we obtain
\begin{align}
    |(IR-R,b_1G_{mn,x})_{\Omega_{f,xy}}| & \leq \|IR-R\|_{\Omega_{f,xy}} \|b_1 G_{mn,x}\|_{\Omega_{f,xy}} \nonumber\\
    & \leq C \|IR-R\|_{L^{\infty}( \Omega_{f,xy})} \sqrt{|\Omega_{f,xy}|} \; \| G_{mn,x}\|_{\Omega_{f,xy}} \nonumber\\
     &\leq C\|IR-R\|_{L^{\infty}( \Omega_{f,xy})} \log(1/\varepsilon^{3/2})\varepsilon^{1/4}\;\| G_{mn}\|_{1,\varepsilon} \nonumber\\
      &\leq \frac{C\sqrt{N \log(1/{\varepsilon})}}{\varepsilon^{-1/4}} \left[\frac{\varepsilon^2 \log^2(1/\varepsilon)}{N^2} +\frac{\varepsilon^{3/2}\log^2(1/\varepsilon^{3/2})}{ N^2} +\frac{\varepsilon \log^2(1/\varepsilon^{3/2})}{N^2} \right] \nonumber\\
       &\leq C \frac{(\log({1/\varepsilon^{3/2}))^3}}{N^{3/2}}.
\end{align}
Consider the second term of the R.H.S. of (\ref{Convection f,xy}) is $|(IF_1-F_1,b_1G_{mn,x})_{\Omega_{f,xy}}|$. By using Theorem \ref{decompositionu}, Theorem \ref{ajeev}, Corollary \ref{A_3} and furthermore, we obtain
\begin{align}
    |(IF_1-F_1,b_1G_{mn,x})_{\Omega_{f,xy}}| & \leq \|IF_1-F_1\|_{\Omega_{f,xy}} \|b_1 G_{mn,x}\|_{\Omega_{f,xy}} \nonumber\\
    & \leq C \|IF_1-F_1\|_{L^{\infty}( \Omega_{f,xy})} \sqrt{|\Omega_{f,xy}|} \; \| G_{mn,x}\|_{\Omega_{f,xy}} \nonumber\\
    &\leq C\|IF_1-F_1\|_{L^{\infty}( \Omega_{f,xy})} \varepsilon^{3/4}\log(1/\varepsilon^{3/2})\;\| G_{mn,x}\|_{\Omega_{f,xy}} \nonumber\\
      &\leq C\|IF_1-F_1\|_{L^{\infty}( \Omega_{f,xy})}\varepsilon^{1/4} \log^{1/2}(1/\varepsilon)\;\| G_{mn}\|_{1,\varepsilon} \nonumber\\
       &\leq C \sqrt{N} \log(1/{\varepsilon})\left[\frac{\varepsilon^2 \log^2(1/\varepsilon)}{N^2\varepsilon^2} +\frac{\varepsilon^{3/2}\log^2(1/\varepsilon^{3/2})}{\varepsilon N^2} +\frac{\varepsilon \log^2(1/\varepsilon^{3/2})}{N^2} \right] \nonumber\\
       &\leq C \frac{(\log({1/\varepsilon^{3/2}))^3}}{N^{3/2}}.
\end{align}
Next consider the third term of the R.H.S. of (\ref{Convection f,xy}) is $|(IF_2-F_2,b_1G_{mn,x})_{\Omega_{f,xy}}|$. By using Theorem \ref{decompositionu}, Theorem \ref{ajeev}, Corollary \ref{B_3} and furthermore, we obtain
\begin{align}
    |(IF_2-F_2,b_1G_{mn,x})_{\Omega_{f,xy}}| & \leq \|IF_2-F_2\|_{\Omega_{f,xy}} \|b_1 G_{mn,x}\|_{\Omega_{f,xy}} \nonumber\\
    & \leq C \|IF_2-F_2\|_{L^{\infty}( \Omega_{f,xy})} \sqrt{|\Omega_{f,xy}|} \; \| G_{mn,x}\|_{\Omega_{f,xy}} \nonumber\\
      &\leq C\|IF_2-F_2\|_{L^{\infty}( \Omega_{f,xy})}\varepsilon^{1/4} \log(1/\varepsilon^{3/2})\;\| G_{mn}\|_{1,\varepsilon} \nonumber\\
      &\leq C\sqrt{N} \log(1/{\varepsilon}) \left[\frac{\varepsilon^2 \log^2(1/\varepsilon)}{N^2} +\frac{\varepsilon^{3/2}\log^2(1/\varepsilon^{3/2})}{\varepsilon^{1/2} N^2} +\frac{\varepsilon \log^2(1/\varepsilon^{3/2})}{\varepsilon N^2} \right] \nonumber\\
       &\leq C \frac{(\log({1/\varepsilon^{3/2}))^3}}{N^{3/2}}.
\end{align}
Next consider the fourth term of the R.H.S. of (\ref{Convection f,xy}) is $|(IF_{12}-F_{12},b_1G_{mn,x})_{\Omega_{f,xy}}|$. By using Theorem \ref{decompositionu}, Theorem \ref{ajeev}, Corollary \ref{A_3}, Corollary \ref{B_3} and furthermore, we obtain
\begin{align}
    |(IF_{12}-F_{12},b_1G_{mn,x})_{\Omega_{f,xy}}| & \leq \|IF_{12}-F_{12}\|_{\Omega_{f,xy}} \|b_1 G_{mn,x}\|_{\Omega_{f,xy}} \nonumber\\
    & \leq C \|IF_{12}-F_{12}\|_{L^{\infty}( \Omega_{f,xy})} \sqrt{|\Omega_{f,xy}|} \; \| G_{mn,x}\|_{\Omega_{f,xy}} \nonumber\\
      &\leq C\|IF_{12}-F_{12}\|_{L^{\infty}( \Omega_{f,xy})} \varepsilon^{1/4}\log(1/\varepsilon^{3/2})\;\| G_{mn}\|_{1,\varepsilon} \nonumber\\
      &\leq C \varepsilon^{1/4} [h_x^2\|F_{12,xx}\|_{L^{\infty}( \Omega_{f,xy})}+h_x h_y\|F_{12,xy}\|_{L^{\infty}( \Omega_{f,xy})}  \nonumber\\
       &  \quad\quad\quad\quad\quad+h_y^2\|F_{12,yy}\|_{L^{\infty}( \Omega_{f,xy})}] (N \log(1/{\varepsilon}))^{1/2}\nonumber\\
     &\leq C\sqrt{N \log(1/{\varepsilon})} \left[\frac{\varepsilon^2 \log^2(1/\varepsilon)}{\varepsilon^2 N^2} +\frac{\varepsilon^{3/2}\log^2(1/\varepsilon^{3/2})}{\varepsilon^{3/2} N^2} +\frac{\varepsilon \log^2(1/\varepsilon^{3/2})}{\varepsilon N^2} \right] \nonumber\\
       &\leq C \left[\frac{ \log^2(1/\varepsilon)}{N^{3/2}} +\frac{\log^2(1/\varepsilon^{3/2})}{ N^{3/2}} +\frac{ \log^2(1/\varepsilon^{3/2})}{ N^{3/2}} \right] 
       \leq C \frac{(\log({1/\varepsilon^{3/2}))^3}}{N^{3/2}}.
\end{align}
Hence, when source point $(x_m,y_n) \in \Omega_{f,x}$, then we have
$$|(Iu-u,b_1G_{mn,x})_{\Omega_{f,xy}}| \leq C \frac{(\log(1/\varepsilon^{3/2}))^3}{N^{3/2}}.$$
Similar arguments we can show that when source point $(x_m,y_n) \in \Omega_{c}$, then we have
$$|(Iu-u,b_1G_{mn,x})_{\Omega_{f,xy}}| \leq C \frac{(\log(1/\varepsilon^{3/2}))^{5/2}}{N}.$$
Hence, when source point $(x_m,y_n) \in \Omega_{c}\cup \Omega_{f,x}$, we have 
$$|(Iu-u,b_1G_{mn,x})_{\Omega_{f,xy}}| \leq C \frac{(\log(1/\varepsilon^{3/2}))^{5/2}}{N}.$$

\begin{lemma}\label{convection term f,c}
    For $\alpha \geq 5/2$, $\beta\geq 2/3$, the following estimates holds:
    \begin{enumerate}
        \item[(i)] When the source point $(x_m,y_n)\in  \Omega_{f,x} \;\text{then}\;|(Iu-u,b_1G_{mn,x})_{\Omega_{c}}| \leq C \frac{(\log(1/\varepsilon))^{1/2}}{N^{1/2}},$
        \item[(ii)] when the source point $\;(x_m,y_n)\in  \Omega_{c} \;\text{then} ~|(Iu-u,b_1G_{mn,x})_{\Omega_{c}}|\leq C.$
    \end{enumerate}
\end{lemma}
\textbf{Proof.} (i) Take the source point $(x_m,y_n)\in \Omega_{f,x}$, consider $ (b_1|(Iu-u)_x,G_{mn})_{\Omega_{c}}| $ and decompose the domain in four parts, we have
\begin{align}
    |(Iu-u,b_1G_{mn,x})_{\Omega_{c}}| \leq & |(IR-R,b_1G_{mn,x})_{\Omega_{c}}| +  |(IF_1-F_1,b_1G_{mn,x})_{\Omega_{c}}|\nonumber\\
    + &  |(IF_2-F_2,b_1G_{mn,x})_{\Omega_{c}}|
    + |(IF_{12}-F_{12},b_1G_{mn,x})_{\Omega_{c}}|. \label{Convection c}
\end{align}
Consider the first term of the R.H.S. of (\ref{Convection c}) is $|(IR-R,b_1G_{mn,x})_{\Omega_{c}}|$. By using inverse inequality, Theorem \ref{decompositionu}, Theorem \ref{ajeev} and furthermore, we obtain
\begin{align}
    |(IR-R,b_1G_{mn,x})_{\Omega_{c}}|& \leq \|IR-R\|_{\Omega_{c}} \|b_1 G_{mn,x}\|_{\Omega_{c}} 
     \leq C \|IR-R\|_{L^{\infty}( \Omega_{c})} \sqrt{|\Omega_{c}|} \; \| G_{mn,x}\|_{\Omega_{c}} \nonumber\\
     &\leq CN\|IR-R\|_{L^{\infty}( \Omega_{c})} \| G_{mn}\|_{\Omega_c} \nonumber\\
      &\leq CN\sqrt{N \log(1/{\varepsilon})} [h_x^2\|R_{xx}\|_{L^{\infty}( \Omega_{c})}+h_x h_y\|R_{xy}\|_{L^{\infty}( \Omega_{c})}+h_y^2\|R_{yy}\|_{L^{\infty}( \Omega_{c})}]  \nonumber\\
      &\leq CN\sqrt{N \log(1/{\varepsilon})}\left[\frac{1}{N^2} +\frac{1}{N^2} +\frac{1}{N^2} \right]  \leq \frac{C( \log(1/{\varepsilon}))^{1/2}}{N^{1/2}}.
\end{align}
Consider the second term of the R.H.S. of (\ref{Convection c}) is $|(IF_1-F_1,b_1G_{mn,x})_{\Omega_{c}}|$. By using Theorem \ref{decompositionu}, Theorem \ref{ajeev}, Lemma \ref{A_1} and furthermore, we obtain
\begin{align}
    |(IF_1-F_1,b_1G_{mn,x})_{\Omega_{c}}|& \leq \|IF_1-F_1\|_{\Omega_{c}} \|b_1 G_{mn,x}\|_{\Omega_{c}}
     \leq C \|IF_1-F_1\|_{L^{\infty}( \Omega_{c})} \sqrt{|\Omega_{c}|} \; \| G_{mn,x}\|_{\Omega_{c}} \nonumber\\
    &\leq C\|IF_1-F_1\|_{L^{\infty}( \Omega_{c})} \;\| G_{mn,x}\|_{\Omega_{c}} \nonumber\\
      &\leq C\|IF_1-F_1\|_{L^{\infty}( \Omega_{c})})\;\varepsilon^{-1/2}\| G_{mn}\|_{1,\varepsilon} \nonumber\\
        &\leq \frac{C\sqrt{N \log(1/{\varepsilon})}}{\sqrt{\varepsilon}} \left[\frac{\|e^{\frac{-\alpha x}{\varepsilon}}\|_{L^{\infty}( \Omega_{c})}}{\varepsilon^2N^2} +\frac{\|e^{\frac{-\alpha x}{\varepsilon}}\|_{L^{\infty}( \Omega_{c})}}{\varepsilon N^2} +\frac{\|e^{\frac{-\alpha x}{\varepsilon}}\|_{L^{\infty}( \Omega_{c})}}{N^2} \right]\nonumber\\
        &\leq  \frac{C\sqrt{ \log(1/{\varepsilon})}}{N^{3/2}\varepsilon^{5/2}} \varepsilon^\alpha \leq \frac{C( \log(1/{\varepsilon}))^{1/2}}{N^{3/2}}\quad \text{when} \; \alpha \geq 5/2.
\end{align}
Next consider the third term of the R.H.S. of (\ref{Convection c}) is $|(IF_2-F_2,b_1G_{mn,x})_{\Omega_{c}}|$. By using inverse inequality, Theorem \ref{decompositionu}, Theorem \ref{ajeev}, Lemma \ref{B_1} and furthermore, we obtain
\begin{align}
    |(IF_2-F_2,b_1G_{mn,x})_{\Omega_{c}}|& \leq \|IF_2-F_2\|_{\Omega_{c}} \|b_1 G_{mn,x}\|_{\Omega_{c}} \nonumber\\
    & \leq C \|IF_2-F_2\|_{L^{\infty}( \Omega_{c})} \sqrt{|\Omega_{c}|} \; \| G_{mn,x}\|_{\Omega_{c}} \nonumber\\
      &\leq C\|IF_2-F_2\|_{L^{\infty}( \Omega_{c})} \;\| G_{mn}\|_{\Omega_c}\nonumber\\
       &\leq CN\sqrt{N \log\frac{1}{\varepsilon}} \left[\frac{1}{N^2} +\frac{1}{\varepsilon^{1/2} N^2} +\frac{1}{\varepsilon N^2} \right] \|e^{-\frac{\beta}{\sqrt{\varepsilon}}(1+y)}+e^{-\frac{\beta}{\sqrt{\varepsilon}}(1-y)}\|_{L^\infty(\Omega_{c})} \nonumber\\
       &\leq  \frac{C}{\varepsilon N^{1/2}} \sqrt{\log(1/{\varepsilon})}\varepsilon^{\frac{3\beta}{2}} \leq \frac{C( \log(1/{\varepsilon}))^{1/2}}{N^{1/2}}\quad \text{when} \; \beta\geq 2/3.
\end{align}
Next consider the fourth term of the R.H.S. of (\ref{Convection c}) is $|(IF_{12}-F_{12},b_1G_{mn,x})_{\Omega_{c}}|$. By using inverse inequality, Theorem \ref{decompositionu}, Theorem \ref{ajeev}, Lemma \ref{A_1}, Lemma \ref{B_1} and furthermore, we obtain
\begin{align}
    |(IF_{12}-F_{12},b_1G_{mn,x})_{\Omega_{c}}|& \leq \|IF_{12}-F_{12}\|_{\Omega_{c}} \|b_1 G_{mn,x}\|_{\Omega_{c}} \nonumber\\
    & \leq C \|IF_{12}-F_{12}\|_{L^{\infty}( \Omega_{c})} \sqrt{|\Omega_{c}|} \; \| G_{mn,x}\|_{\Omega_{c}} \nonumber\\
      &\leq CN\|IF_{12}-F_{12}\|_{L^{\infty}( \Omega_{c})} \;\| G_{mn}\|_{\Omega_c} \nonumber\\
    &\leq CN\sqrt{N \log\frac{1}{\varepsilon}}  \left[\frac{\varepsilon^{-2}}{N^2} +\frac{\varepsilon^{-3/2}}{ N^2} +\frac{1}{\varepsilon N^2} \right]\varepsilon^{\alpha+\frac{3\beta}{2}}\nonumber\\
       &\leq C \varepsilon^{\alpha+\frac{3\beta}{2}-2}\frac{(\log(1/\varepsilon))^{1/2}}{N^{1/2}}  \leq  \frac{C(\log(1/\varepsilon))^{1/2}}{N^{1/2}} \quad \text{when} \; \alpha+\frac{3\beta}{2}\geq2.
\end{align}
Hence, when source point $(x_m,y_n) \in \Omega_{f,x}$, we have
$$ |(Iu-u,b_1G_{mn,x})|  \leq C \frac{(\log(1/\varepsilon))^{1/2}}{N^{1/2}}.$$
(ii) Similarly when the source point $(x_m,y_n) \in \Omega_{c} $, using Lemma \ref{Greens in c} yields $\|G_{mn}\|\leq CN$, we have 
\begin{center}
 $|(Iu-u,b_{1}G_{mn,x})|  \leq C.$
 \end{center}
\begin{Remark}
     When the source point $(x_m,y_n)\in \Omega_c$, we can obtain only $\varepsilon$-independent bounds; the present approach does not provide a convergence rate as reported \cite{Stynes}.
\end{Remark}
Following a similar argument to that discussed in the previous lemmas (\ref{convection term f,x} - \ref{convection term f,c}), we obtained the next Lemma. 
\begin{lemma}\label{reaction full}
    For $\alpha \geq 5/2$, $\beta\geq 2/3$, the following estimates holds:
    \begin{enumerate}
        \item[(i)] When the source point $(x_m,y_n)\in  \Omega_{f,x} \;\text{then}\;|((c-b_{1,x}),(Iu-u)G_{mn})| \leq C \frac{\log(1/\varepsilon)}{N^{3/2}}$,
        \item[(ii)] when the source point $\;(x_m,y_n)\in  \Omega_{c} \;\text{then} \; |((c-b_{1,x}),(Iu-u)G_{mn})| \leq C \frac{(\log(1/\varepsilon^{3/2}))^{7/2}}{N}.$
    \end{enumerate}
\end{lemma}
 \section{Main Result}\label{section6}
 Using the analysis from the preceding sections, we derive our main result, which confirms an $O(N^{-1/2}(\log(1/\varepsilon))^{1/2})$ rate of convergence within the $x$-layer region. For the coarse layer region, the result is $\varepsilon$-independent, but we could not determine a specific order of convergence. 
\begin{theorem} \label{main theorem}
   Let  $u_h$ denote the approximate solution of the problem \ref{current problem} and  the source point $(x_m,y_n) \in \Omega_{f,x}$, then for $\alpha\geq 5/2$ and $\beta\in[2/3,1]$ the following estimate hold
\begin{equation*}
    |(u-u_h)(x_m,y_n)|\leq C \frac{(\log(1/\varepsilon))^{1/2}}{N^{1/2}}.
\end{equation*}
\end{theorem}
\textbf{Proof.}
Let $(x_m,y_n)$ be a mesh point, then
\begin{equation}
    (Iu-u_h)(x_m,y_n) = (\delta_{mn},Iu-u_h),
\end{equation}
where $\delta_{mn}$ is the Kronecker delta function defined on mesh point $(x_m,y_n)$, then
\begin{align*}
 (Iu-u_h)(x_m,y_n) &= (\delta_{mn},Iu-u_h) = B(Iu-u_h,G_{mn}),\\
  \text{where} \quad B(u,v) &= \varepsilon \int_\Omega \nabla u.\nabla v \;dx + \int_\Omega b_1 u_x v \;dx + \int_\Omega c uv \; dx\\
    &= \varepsilon ( \nabla u,\nabla v ) + ( b_1 u_x, v) + (c u,v ).\\
    \implies  (Iu-u_h)(x_m,y_n) &= B(Iu-u_h,G_{mn}) \\
    & =B(Iu-u,G_{mn}) + B(u-u_h,G_{mn}).
\end{align*}
But Galerkin's orthogonality implies $B(u-u_h,G_{mn}) = 0, $ therefore 
\begin{align*}
&(Iu-u_h)(x_m,y_n) = B(Iu-u,G_{mn})\\
    &  = \varepsilon \int_\Omega \nabla(Iu-u).\nabla G_{mn} \;dx dy + \int_\Omega b_1 (Iu-u)_x G_{mn} \;dx dy  + \int_\Omega c(Iu-u)G_{mn} \;dx dy\\
    &= \varepsilon  (\nabla(Iu-u),\nabla G_{mn}) + ( b_1 (Iu-u)_x, G_{mn}) + (c(Iu-u),G_{mn}).
\end{align*}
Applying integration by parts yields
\begin{align*} 
      (Iu-u_h)(x_m,y_n) = \varepsilon  (\nabla(Iu-u),\nabla G_{mn}) - ( b_1 (Iu-u), G_{mn,x})+ ((c-b_{1,x}),(Iu-u)G_{mn})
\end{align*} 
therefore,
\begin{align}\label{Eq:8}
    |(Iu-u_h)(x_m,y_n)| \;\leq |\varepsilon  (\nabla(Iu-u),\nabla G_{mn})| + |( b_1 (Iu-u), G_{mn,x})| + |((c-b_{1,x}),(Iu-u)G_{mn})|. 
\end{align}
To obtain the error estimate, we aim to estimate the term on the R.H.S.
of (\ref{Eq:8}). Consider the diffusion term of (\ref{Eq:8}) and decompose the solution $u$ into regular and singular parts as described in Theorem \ref{decompositionu}. We obtain 
\begin{multline} \label{eq:10}
    |\varepsilon(\nabla(u-Iu),\nabla G_{mn})|\;  = |\varepsilon(\nabla((R +F_1 + F_2 +F_{12})-I(R +F_1 + F_2 +F_{12})),\nabla G_{mn})|\\
    \leq|\varepsilon (\nabla(R-IR),\nabla G_{mn})|+|\varepsilon(\nabla(F_1-IF_1),\nabla G_{mn})|\\
    +|\varepsilon(\nabla(F_2-IF_2),\nabla G_{mn})|+|\varepsilon(\nabla(F_{12}-IF_{12}),\nabla G_{mn})| 
\end{multline}
Using Lemmas (\ref{diffusion R} - \ref{diffusion F_12}), for the source point $(x_m,y_n) \in \Omega_{c} \cup\Omega_{f,x}$, then the diffusion term is bounded by
\begin{equation}\label{final diffusion}
    |(\nabla(u-Iu),   \nabla G_{mn})| \leq C \frac{(\log(1/\varepsilon^{3/2})^{7/2}}{N}.
\end{equation} 
Consider the second term $|( b_1 (Iu-u), G_{mn,x})|$ of (\ref{Eq:8}). Applying Lemmas (\ref{convection term f,x} - \ref{convection term f,c}) for any source point $(x_m,y_n) \in \Omega_{f,x}$, we obtain
\begin{align}\label{final convection}
    |(Iu-u,b_1G_{mn,x})|  \leq C \frac{(\log(1/\varepsilon))^{1/2}}{N^{1/2}}.
\end{align}
Next consider $|((c-b_{1,x}),(Iu-u)G_{mn})|$, i.e. the third term of (\ref{Eq:8}), Lemma \ref{reaction full} yields for any     source point $(x_m,y_n) \in \Omega_{f,x}$, the following estimate hold
\begin{align}\label{final reaction}
    |((c-b_{1,x}),(Iu-u)G_{mn})|  \leq C \frac{\log(1/\varepsilon)}{N^{3/2}}.
\end{align}
Combining \eqref{final diffusion}, \eqref{final convection}, and \eqref{final reaction} yields the result.
\begin{Remark}
    Applying Lemmas (\ref{diffusion R} - \ref{diffusion F_12}), Lemmas (\ref{convection term f,x} - \ref{convection term f,c}), Lemma \ref{reaction full} and Green's function estimate Lemma \ref{Greens in c} for the source point $(x_m,y_n)\in \Omega_c$, we have $|(u-u_h)(x_m,y_n)|\leq C.$ Here we obtain an $\varepsilon$-independent bound for the estimate; however no convergence is obtained and it agrees with the findings of \cite{Stynes}.
\end{Remark}
\begin{Remark}
     For the coarse layer region, we have not yet found an order of convergence, but we believe higher-degree polynomials will help us to obtain it.  We explore it further in future work.
\end{Remark}
\section{Numerical results and discussion}\label{section-6}
In this section, we take a numerical example to validate the theoretical results established in Theorems (\ref{main theorem}). We apply our method to a test problem to verify the uniform rate of convergence only in the $x$-layer and coarse layer regions.
\subsection{Numerical implementation}
We evaluate the error and the uniform rate of convergence with respect to the maximum norm in the coarse layer and $x$-layer regions for the different values of the perturbation parameters $\varepsilon$ and the discretization parameter N. In the absence of an exact solution to the text example, we employ the double mesh principle technique to compute the error and order of the convergence. We calculate the error estimates in the maximum norm using $\| U_N-U_{2N}\|_{L^\infty{(\Omega})}$. The convergence rate is calculated by using the following:
$$C^N_\varepsilon = \log_2\left(\frac{\|U_N-U_{2N}\|_{L^\infty(\Omega)}}{\|U_
{2N}-U_{4N}\|_{L^\infty(\Omega)}}\right).$$
A numerical example is presented where the exact solution is not available.\\
\begin{example}\label{example}
        We investigate the turning point singularly perturbed equation provided below
        \begin{equation*}
\begin{split}
-\varepsilon \Delta u - x (x^2 +e^{1+y^2}) u_x + (2- x^2 - y^2) u &= f(x,y), \quad \text{in} \; \Omega = (-1,1)\times(-1,1) \\
u &=0, \quad \Gamma = \partial \Omega. 
\end{split}
\end{equation*}
\end{example}
For this particular problem, the source term is set as $f(x,y)=\frac{xy}{1+x^2+y^2}$.
A key feature of this problem is its turning point at $x=0$, which gives rise to an interior layer along the $x$-axis. In parallel, two boundary layers are present at the bottom of the boundary $y=-1$ and the top of the boundary $y=1$ of the domain. A double-mesh principle is used to improve the accuracy and convergence behavior of the numerical scheme. This scheme allows us to estimate errors and evaluate the rate of convergence for different values of the perturbation parameter $\varepsilon$ and the mesh discretization parameter $N$ without requiring the exact solution. The resultant error values and computed rates of convergence are systematically presented in tables.
\begin{table}[ht]
\centering
\caption{Error in the subdomain $\Omega_c$ for example \eqref{example}.}
\label{tab:table1}
\begin{tabular}{|c|ccccc|}
\hline
 & \multicolumn{5}{c|}{Error in $\Omega_c$}  \\\cline{2-6}  
 & \multicolumn{5}{c|}{Number of discretization parameter $N$} \\ \hline
 $\varepsilon$  & 32 & 64&128 & 256& 512\\ \hline
    $10^{-5}$  & 9.113E-03 & 1.062E-02 & 1.077E-02 & 1.020E-02 & 1.190E-03 \\
        $10^{-6}$  & 8.664E-03 & 1.039E-02 & 1.042E-02 & 8.919E-03 & 7.338E-03 \\
        $10^{-7}$  & 8.148E-03 & 1.043E-02 & 1.084E-02 & 9.427E-03 &7.093E-03\\
        $10^{-8}$  & 7.794E-03 & 1.043E-02 & 1.127E-02 & 1.018E-02 &7.748E-03  \\
        $10^{-9}$  & 7.380E-03 & 1.032E-02 & 1.159E-02 & 1.084E-02&8.467E-03  \\\hline
\end{tabular}
\end{table}
\begin{table}
\centering
\caption{Error in the subdomain $\Omega_{f,x}$ for example \eqref{example}.}
\label{tab:table2}
\begin{tabular}{|c|ccccc|}
\hline
 & \multicolumn{5}{c|}{Error in $\Omega_{f,x}$}  \\\cline{2-6}  
 & \multicolumn{5}{c|}{Number of discretization parameter $N$} \\ \hline
 $\varepsilon$  & 32 & 64&128 & 256&512\\ \hline
     $10^{-5}$  & 1.878E-02 & 1.494E-02 & 1.199E-02 & 1.102E-02&7.665E-03 \\
        $10^{-6}$ & 1.951E-02 & 1.647E-02 & 1.141E-02 & 9.324E-03& 7.565E-03 \\
        $10^{-7}$  & 2.001E-02 & 1.811E-02 & 1.234E-02 & 9.803E-03 & 7.190E-03\\
        $10^{-8}$  & 2.028E-02 & 1.750E-02 & 1.383E-02 & 1.089E-02 & 7.812E-03 \\
        $10^{-9}$  & 2.039E-02 & 1.903E-02 & 1.526E-02 & 1.183E-02 & 8.547E-03 \\\hline
\end{tabular}
\end{table}
\begin{figure}
    \centering
    \includegraphics[width=0.5\linewidth]{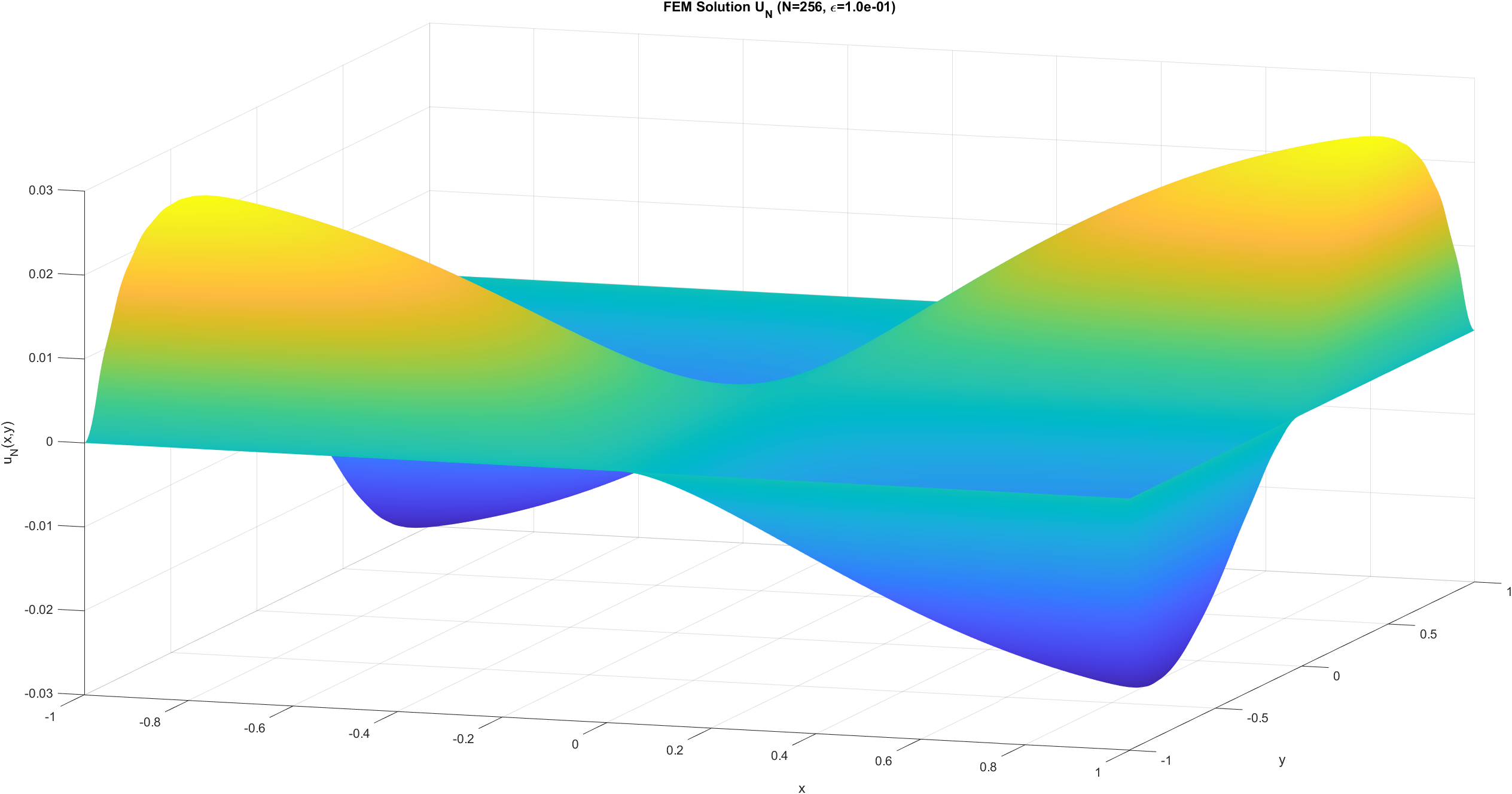}
    \label{solution plot}
    \caption{FEM solution $U_N$(N = 256, $\varepsilon$ = 1.0e-01)}
\end{figure}
\begin{table}
\centering
\caption{Rate of convergence in the subdomains $\Omega_c$ and $\Omega_{f,x}$ for test example \eqref{example}.}
\label{tab:table3}
\begin{tabular}{|c|cccc|cccc|}
\hline
 & \multicolumn{4}{c|}{$\Omega_c$} & \multicolumn{4}{c|}{$\Omega_{f,x}$} \\ \cline{2-9} 
 & \multicolumn{4}{c|}{Number of discretization parameter $N$} & \multicolumn{4}{c|}{Number of discretization parameter $N$} \\ \hline
$\varepsilon$ & 32 & 64 & 128&256 & 32 & 64 & 128 &256\\ \hline
         $10^{-5}$  & -0.2213 & -0.0197 & 0.0788 & 0.2703  & 0.3294 & 0.3178 & 0.1208 & 0.5244 \\
        $10^{-6}$  & -0.2623 & -0.0043 & 0.2246 & 0.2815  & 0.2440 & 0.5303 & 0.2907 & 0.3015 \\
        $10^{-7}$  & -0.3563 & -0.0557 & 0.2017 & 0.4104  & 0.1437 & 0.5529 & 0.3327 & 0.4472 \\
        $10^{-8}$  & -0.4204 & -0.1111 & 0.1468 & 0.3931  & 0.1567 & 0.4951 & 0.3456 & 0.4788 \\
        $10^{-9}$  & -0.4833 & -0.1682 & 0.0962 & 0.3571  & 0.1165 & 0.4352 & 0.3668 & 0.4690 \\\hline
\end{tabular}
\end{table}
\begin{table}
\centering
\caption{Error in the subdomain $\Omega_{f,y}$ for example \eqref{example}.}
\label{tab:table6}
\begin{tabular}{|c|ccccc|}
\hline
 & \multicolumn{5}{c|}{Error in $\Omega_{f,y}$}  \\\cline{2-6}  
 & \multicolumn{5}{c|}{Number of discretization parameter $N$} \\ \hline
 $\varepsilon$ & 32 & 64&128 & 256&512\\ \hline
    $10^{-5}$  & 1.090E-01 & 1.276E-01 & 1.456E-01 & 1.462E-01&8.547E-03 \\
        $10^{-6}$  & 1.089E-01 & 1.286E-01 & 1.526E-01 & 1.794E-01& 1.990E-01\\
        $10^{-7}$  & 1.089E-01 & 1.284E-01 & 1.530E-01 & 1.832E-01&2.186E-01 \\
        $10^{-8}$  & 1.088E-01 & 1.281E-01 & 1.528E-01 & 1.833E-01 &2.206E-01 \\
        $10^{-9}$  & 1.086E-01 & 1.278E-01 & 1.526E-01 & 1.831E-01& 2.206E-01 \\\hline
\end{tabular}
\end{table}
\begin{table}
\centering
\caption{Error in the subdomain $\Omega_{f,xy}$ for example \eqref{example}.}
\label{tab:table4}
\begin{tabular}{|c|ccccc|}
\hline
 & \multicolumn{5}{c|}{Error in $\Omega_{f,xy}$}  \\\cline{2-6}  
 & \multicolumn{5}{c|}{Number of discretization parameter $N$} \\ \hline
 $\varepsilon$ & 32 & 64&128 & 256&512\\ \hline
   $10^{-5}$  & 1.084E-01 & 1.274E-01 & 1.447E-01 & 1.469E-01& 9.961E-02\\
        $10^{-6}$  & 1.071E-01 & 1.288E-01 & 1.457E-01 & 1.751E-01&1.983E-01 \\
        $10^{-7}$  & 1.063E-01 & 1.289E-01 & 1.454E-01 & 1.751E-01 &2.130E-01\\
        $10^{-8}$  & 1.058E-01 & 1.287E-01 & 1.494E-01 & 1.740E-01 &2.132E-01\\
        $10^{-9}$ & 1.054E-01 & 1.285E-01 & 1.520E-01 & 1.756E-01& 2.126E-01\\\hline
\end{tabular}
\end{table}
\begin{table}
\centering
\caption{Rate of convergence in the subdomains $\Omega_y$ and $\Omega_{f,xy}$ for test example \eqref{example}.}
\label{tab:table5}
\begin{tabular}{|c|cccc|cccc|}
\hline
 & \multicolumn{4}{c|}{$\Omega_y$} & \multicolumn{4}{c|}{$\Omega_{f,xy}$} \\ \cline{2-9} 
 & \multicolumn{4}{c|}{Number of discretization parameter $N$} & \multicolumn{4}{c|}{Number of discretization parameter $N$} \\ \hline
$\varepsilon$ & 32 & 64 & 128&256 & 32 & 64 & 128 &256\\ \hline
$10^{-5}$ & -0.2269 & -0.1905 & -0.0060 & 0.7163 & -0.2340 & -0.1833 & -0.0222 & 0.5609 \\
        $10^{-6}$ & -0.2403 & -0.2459 & -0.2335 & -0.1502 & -0.2663 & -0.1781 & -0.2650 & -0.1794 \\
        $10^{-7}$ & -0.2375 & -0.2529 & -0.2592 & -0.2550 & -0.2780 & -0.1738 & -0.2687 & -0.2822 \\
        $10^{-8}$ & -0.2356 & -0.2549 & -0.2627 & -0.2667 & -0.2831 & -0.2152 & -0.2195 & -0.2934 \\
        $10^{-9}$ & -0.2345 & -0.2559 & -0.2631 & -0.2687 & -0.2856 & -0.2420 & -0.2089 & -0.2758 \\\hline
\end{tabular}
\end{table}
\subsection{Validation} For text example \eqref{example}, Tables \eqref{tab:table1}, \eqref{tab:table2}, \eqref{tab:table6} and \eqref{tab:table4} display the numerical errors for subdomains $\Omega_{c}$, $\Omega_{f,x}$, $\Omega_{f,y}$ and $\Omega_{f,xy}$ respectively. Additionally, Tables \eqref{tab:table3} and \eqref{tab:table5} present the rates of convergence when the source point is located within $\Omega_{c}\cup\Omega_{f,x}$ and $\Omega_{f,y}\cup\Omega_{f,xy}$, respectively. As shown in Table \eqref{tab:table3}, the proposed scheme is uniformly convergent for the source point that lies in  $\Omega_{f,x}$. Conversely, Table \eqref{tab:table5} indicates that uniformity is lost when the source point lies in $\Omega_{f,y}\cup\Omega_{f,xy}$. The solution plot clearly captures the interior and boundary layers. Furthermore, we obtain the numerical convergence rates, which are $N^{-1/2}$ up to a logarithmic factor when the source point lies within the subdomain $\Omega_{f,x}$, confirming our theoretical predictions.
\section{Conclusions}
In this article, we employ a standard finite element method to solve a two-dimensional singularly perturbed convection-diffusion problem with a turning point. Because the convection coefficient changes sign within the domain, the problem has an interior layer parallel to the $y$-axis at $x=0$ and two boundary layers parallel to the $x$-axis at $y=-1$ and  $y=1$. Parameter uniform error estimates are obtained for $x$-layer regions. For the coarse region, we obtain only $\varepsilon$-independent estimates, similar to those reported \cite{Stynes}, which uses a Shishkin-type mesh with a specific mesh transition parameter. Additionally, pointwise error estimates are established through the application of the discrete Green’s function. We prove uniform convergence in the perturbation parameter of order $N^{-1/2}(\log(1/\varepsilon))^{1/2}$ pointwise in the $x$-layer region.
\section*{Acknowledgments}
The authors are thankful to the anonymous referees for their valuable comments and suggestions to improve the paper. The authors thank LNMIIT, Jaipur, for providing the facilities and support necessary to complete this research. The authors also acknowledge the DST-FIST program (Govt. of India) for providing the financial support for setting up the computing lab facility under the scheme "Fund for Improvement of Science and Technology" (FIST - No. SR/EST/MS-1/2018/24). The authors also acknowledge the technical support provided by the LNM Institute of Information Technology, India.


\begin{thebibliography}{99}
\bibitem{Aasna2024}
Aasna Rai, P.: SDFEM for time-dependent singularly perturbed boundary turning point problem. Int. J. Comput. Methods (2024).
\bibitem{Aasna Rai2}
Aasna Rai, P.: Analysis of SDFEM for singularly perturbed delay differential equation with boundary turning point. Int. J. Appl. Comput. Math. 10, 20 (2024).
\bibitem{Al-Haboobi_24}
Al-Haboobi, A., Al-Juaifri, G., Al-Muslimawi, A.: Numerical study of Newtonian laminar Flow around circular and square cylinders. Results Control Optim. 14, 100328 (2024).
\bibitem{Becher_25}
Becher, S., Roos, H.: Richardson extrapolation for a singularly perturbed turning point problem with exponential boundary layers. J. Comput. Appl. Math. 290, 334–351 (2015).
\bibitem{Becher_16}
Becher, S.: FEM-analysis on graded meshes for turning point problems exhibiting an interior layer. ArXiv:1603.04653 (2016).
\bibitem{Becher18}
Becher, S.: Analysis of Galerkin and streamline-diffusion FEMs on piecewise equidistant meshes for turning point problems exhibiting an interior layer. Appl. Numer. Math. 123, 121–136 (2018).
\bibitem{Berger_84}
Berger, A., Han, H., Kellogg, R.: A priori estimates and analysis of a numerical method for a turning point problem. Math. Comput. 42, 465–492 (1984).
\bibitem{Chen08}
Chen, L., Wang, Y., Wu, J.: Stability of a streamline diffusion finite element method for turning point problems. J. Comput. Appl. Math. 220, 712–724 (2008).
\bibitem{Farrell_88}
Farrell, P.: Sufficient conditions for the uniform convergence of a difference scheme for a singularly perturbed turning point problem. SIAM J. Numer. Anal. 25, 618–643 (1988).
\bibitem{Feistauer_03}
Feistauer, M., Felcman, J., Straškraba, I.: Mathematical and Computational Methods for Compressible Flow. Oxford University Press, Oxford (2003).
\bibitem{Geng_13}
Geng, F., Qian, S.: Reproducing kernel method for singularly perturbed turning point problems having twin boundary layers. Appl. Math. Lett. 26, 998–1004 (2013).
\bibitem{Gupta_21}
Gupta, V., Sahoo, S., Dubey, R.: Robust higher order finite difference scheme for singularly perturbed turning point problem with two outflow boundary layers. Comput. Appl. Math. 40, 1–23 (2021).
\bibitem{Hahn_87}
Hahn, S., Bigeon, J., Sabonnadiere, J.: An ‘upwind’finite element method for electromagnetic field problems in moving media. Int. J. Numer. Methods Eng. 24, 2071–2086 (1987).
\bibitem{Kadalbajoo_10}
Kadalbajoo, M., Gupta, V.: A parameter uniform B-spline collocation method for solving singularly perturbed turning point problem having twin boundary layers. Int. J. Comput. Math. 87, 3218–3235 (2010).
\bibitem{Kumar_19}
Kumar, D.: A parameter-uniform method for singularly perturbed turning point problems exhibiting interior or twin boundary layers. Int. J. Comput. Math. 96, 865–882 (2019).
\bibitem{lins}
Linß, Torsten.: Layer-adapted meshes for reaction-convection-diffusion problems, Lecture Notes in Mathematics, Vol. 1985, Springer-Verlag, Berlin, Heidelberg, 2009.
\bibitem{Liseikin18}
Liseikin, V.: Layer Resolving Grids and Transformations for Singular Perturbation Problems. Walter de Gruyter GmbH and Co KG, De Gruyter (2018).
\bibitem{Mbayi22}
Mbayi, Charles K., Justin B. Munyakazi, and Kailash C. Patidar. "A fitted numerical method for interior-layer parabolic convection–diffusion problems." International Journal of Computational Methods 19.10 (2022): 2250028.
\bibitem{Natesan_98}
Natesan, S., Ramanujam, N.: A computational method for solving singularly perturbed turning point problems exhibiting twin boundary layers. Appl. Math. Comput. 93, 259–275 (1998).
\bibitem{Natesan_03}
Natesan, S., Jayakumar, J., Vigo-Aguiar, J.: Parameter uniform numerical method for singularly perturbed turning point problems exhibiting boundary layers. J. Comput. Appl. Math. 158, 121–134 (2003).
\bibitem{Rajeev_Ranjan}
Ranjan, K.R., Gowrisankar, S. Error analysis of discontinuous Galerkin methods on layer-adapted meshes for the two-dimensional turning point problem. J. Appl. Math. Comput. 70, 2453–2485 (2024).
\bibitem{Tobiska_book}
Roos, H.G., Stynes, M., Tobiska, L.: Robust numerical methods for singularly perturbed differential equations: convection-diffusion-reaction and flow problems. Springer, New York (2008).
\bibitem{Stynes}
Stynes, Martin, and Eugene O'Riordan. "A uniformly convergent Galerkin method on a Shishkin mesh for a convection-diffusion problem." Journal of Mathematical Analysis and Applications 214.1 (1997): 36-54.
\bibitem{Sun94}
Sun, G., Stynes, M.: Finite element methods on piecewise equidistant meshes for interior turning point problems. Numer. Algorithms. 8, 111–129 (1994).
\bibitem{Yasir_20}
Yasir, R., Al-Muslimawi, A., Jassim, B.: Numerical simulation of non-Newtonian inelastic flows in channel based on artificial compressibility method. J. Appl. Comput. Mech. 6, 271–283 (2020).
\bibitem{Zarin}
Zarin, H., Roos, H.G.: Interior penalty discontinuous approximations of convection-diffusion problems
with parabolic layers. Numer. Math. 100(4), 735–759 (2005).
\bibitem{zhang}
Zhang Z.: Pointwise error estimates of streamline diffusion finite element methods Numerische Mathematik, (94) (2003), pp. 555-579.
 \end{thebibliography}
\end{document}